\documentclass{article}

\usepackage{indentfirst}
\usepackage{amssymb}
\usepackage{amsmath}
\usepackage{amsfonts}
\usepackage{theorem}

\makeatletter
\def\rightharpoonupfill@{\arrowfill@\relbar\relbar\rightharpoonup}
\newcommand{\xrightharpoonup}[2][]{\ext@arrow 0359\rightharpoonupfill@{#1}{#2}}
\makeatother

\newcommand{\ds}{\displaystyle}

\newtheorem{theo}{Theorem}[section]
\newtheorem{lem}[theo]{Lemma}
\newtheorem{prop}[theo]{Proposition}

\theorembodyfont{\upshape}
\newtheorem{rem}[theo]{Remark}

\numberwithin{equation}{section}

\begin{document}
\title{
\textsc{Spatial heterogeneity in 3D-2D dimensional reduction}}
\author{}
\date{}
\maketitle
\begin{center}
{\sc{Jean-Fran\c{c}ois Babadjian(*),  Gilles A. Francfort(**)}}
\end{center}

\begin{small}
\textbf{Abstract. }A justification of heterogeneous membrane models as zero-thickness limits of a cylindral three-dimensional heterogeneous nonlinear hyperelastic body is proposed in the spirit of \cite{Le&Ra}. Specific characterizations of the 2D elastic energy are produced. As a generalization of \cite{Bo&Fo&Ma}, the case where  external loads induce a density of bending moment that produces a Cosserat vector field is also investigated. Throughout, the 3D-2D dimensional reduction is viewed as a problem of $\Gamma$-convergence of the elastic energy, as the thickness tends to zero.\\

\end{small}

\section{Introduction}

The purpose of this article is to study the behavior of a thin elastic plate, as the thickness tends to zero. This approach renders more realistic the idealized view of a film as a thin plate. The originality of the work comes from the heterogeneity of the material under
consideration. Previous results have been established in the homogeneous case ; our aim here is to generalize those. As we will see, accounting for  inhomogeneity  leads to technical difficulties linked to the equi-integrable character of the scaled gradient. We will use a ``classical'' approach of the theory of dimension reduction. In recent years, the investigation of dimensional reduction has focussed on variational methods and used De Giorgi's $\Gamma$-convergence (see \cite{Bra&De} and \cite{DalMaso}) as its main tool.

As far as 3D-2D asymptotic analysis is concerned, the seminal paper is \cite{Le&Ra}, in which a membrane model is derived
from three-dimensional hyperelasticity. In its footstep several studies have derived or re-derived various membrane-like models in various settings; see in particular \cite{Bra&Fo&Fr}  and references therein; note that in Section 3 of that paper, a transversally inhomogeneous
thin domain is studied, but that in-plane-homogeneity is imposed. Because of frame indifference, it may occur that the membrane effect is not excited by the loads : this is the case for example when the lateral boundary conditions on the thin domain are compressive (see e.g. Theorem 6.2 in \cite{Fr&Ja&Mu2}). Then the membrane energy, which results from a 3D-energy of the order of the thickness  $\varepsilon$,  is actually zero and lower energy modes are activated. In \cite{Fox&Ra&Si}, a justification of classical nonlinear plate models for a homogeneous isotropic material is given by a formal asymptotic expension. Recently, those results have been rigorously justified by means of variational methods for general homogeneous hyperelastic bodies. A Kirchhoff bending model in \cite{Fr&Ja&Mu1} and \cite{Fr&Ja&Mu2}, and a F\"oppl-von K\'arm\'an model in \cite{Fr&Ja&Mu3} have been obtained when the 3D-energy  scales respectively like $\varepsilon^3$ and $\varepsilon^5$.

The present study falls squarely  within the membrane framework in the sense that, thanks to frame indifference, the stored energy function depends only on the first fundamental form of the deformed plate mid-surface. Our goal is to rigorously derive models for heterogeneous membranes from their heterogeneous thin 3D-counterparts. The paper is devoted to a generalization of the results established in \cite{Le&Ra}, \cite{Bra&Fo&Fr} and \cite{Bo&Fo&Ma} to the case of a general inhomogeneity.

The key ingredient of this study is the equi-integrability theorem of \cite{Bo&Fo} (Theorem 1.1 of that reference). An alternative proof
of that theorem  was also communicated to the authors \cite{Braides}. This theorem shows that a sequence of scaled gradients $\left\{\left( D_{\alpha}u_{\varepsilon}|\frac{1}{\varepsilon}D_3u_{\varepsilon}\right)\right\}$, which is bounded in $L^p(\Omega;\mathbb{M}^{3 \times 3})$, with $p>1$, can be decomposed into the sum of two sequences $\{w_{\varepsilon}\}$ and $\{z_{\varepsilon}\}$ where  $\left\{\left| \left(D_{\alpha}w_{\varepsilon}|\frac{1}{\varepsilon}D_3w_{\varepsilon}\right)\right|^p \right\}$   is equi-integrable and $z_{\varepsilon} \to 0$ in measure.

Let $\omega$ be a bounded open subset of $\mathbb{R}^2$. Consider $\Omega_{\varepsilon}:= \omega \times (-\varepsilon,\varepsilon)$, the reference configuration of a hyperelastic  heterogeneous thin film, with elastic energy density  given by the $\varepsilon$-dependent Carath\'eodory function $W_{\varepsilon}:\Omega_{\varepsilon} \times \mathbb{M}^{3 \times 3} \rightarrow \mathbb{R}$. We will assume e.g. that the body is clamped on the lateral boundary  $\Gamma_{\varepsilon}:=\partial\omega\times (-\varepsilon,\varepsilon)$ and that it is submitted to the action of surface traction densities on $\Sigma_\varepsilon:=\omega\times \{-\varepsilon,\varepsilon\}$. The total energy of the system under a deformation $u: \Omega_{\varepsilon} \rightarrow \mathbb{R}^3$ is given by
$$\mathcal{E}(\varepsilon)(u)=\int_{\Omega_{\varepsilon}}W_{\varepsilon}(x;D u) dx - \int_{\Omega_{\varepsilon}}f_{\varepsilon}.u \, dx -\int_{\Sigma_\varepsilon}g_{\varepsilon}.u \, d\mathcal H^2,$$
where $\mathcal H^2$ stands for the 2-dimensional surface measure, $f_{\varepsilon} \in L^{p'}(\Omega_{\varepsilon};\mathbb{R}^3)$ denotes an appropriate dead load and $g_{\varepsilon} \in  L^{p'}(\Sigma_\varepsilon;\mathbb{R}^3)$ some surface traction densities  ($1/p+1/p'=1$). We denote by $W^{1,p}_{\Gamma_{\varepsilon}}(\Omega_{\varepsilon};\mathbb{R}^3)$  the space of kinematically admissible fields, that is the functions in  $W^{1,p}(\Omega_{\varepsilon};\mathbb{R}^3)$ with zero trace on $\Gamma_{\varepsilon}$. As is classical in hyperelasticity, the equilibrium problem is viewed as the minimisation problem
$$\inf_{u-x \in W^{1,p}_{\Gamma_{\varepsilon}}(\Omega_{\varepsilon};\mathbb{R}^3)} \mathcal{E}({\varepsilon})(u).$$
Since the integration domain depends on $\varepsilon$, we  reformulate the problem on a fixed domain through a $1/\varepsilon$-dilatation in the transverse direction $x_3$. Let $x_\alpha$ the vector $(x_1,x_2) \in \omega$, we set $v(x_{\alpha},\frac{x_3}{\varepsilon}):=u(x_{\alpha},x_3)$ and $\mathcal{E}_{\varepsilon}(v):=\frac{1}{\varepsilon}\mathcal{E}(\varepsilon)(u)$, then
\begin{eqnarray*}
\mathcal{E}_{\varepsilon}(v)=\int_{\Omega} W_{\varepsilon}\left( x_{\varepsilon};D_{\alpha} v(x) \Big|\frac{1}{\varepsilon} D_3 v(x) \right) dx - \int_{\Omega} f_{\varepsilon}(x_{\varepsilon}).v(x) dx -\frac{1}{\varepsilon}\int_{\Sigma}g_{\varepsilon}(x_{\varepsilon}).v(x) d\mathcal H^2,
\end{eqnarray*}
where $x_{\varepsilon}:=(x_\alpha,\varepsilon x_3)$. We set  $\Omega:=\omega \times (-1,1)$, $\Sigma:=\omega \times \{-1,1\}$, denote by $D_{\alpha}v$  the $3 \times 2$  matrix of partial derivatives $\frac{\partial v_i}{\partial x_{\alpha}}$ ($i \in \{1,2,3\}$, $\alpha \in \{1,2\}$) and by $(\overline{F}|z)$, the two first columns of which are those of the matrix $\overline{F} \in \mathbb{M}^{3 \times 2}$, while the last one is the vector $z \in \mathbb{R}^3$. A formal asymptotic expension in \cite{Fox&Ra&Si} shows that the membrane theory arises if the body forces is of order 1 and the surfaces loadings is of order $\varepsilon$. We next assume that
$$\left\{
\begin{array}{rcl}
W_{\varepsilon}(x_{\alpha},\varepsilon x_3;F) & = & W(x_{\alpha},x_3;F),\\
f_{\varepsilon}(x_{\alpha},\varepsilon x_3) & = & f(x_{\alpha},x_3),\\
g_{\varepsilon}(x_{\alpha},\varepsilon x_3) & = &  g_0(x_{\alpha},x_3)+\varepsilon g(x_{\alpha},x_3)
\end{array}
\right.$$
where $f \in L^{p'}(\Omega;\mathbb{R}^3)$, $g_0,\, g \in L^{p'}(\Sigma;\mathbb{R}^3)$ and  $W:\Omega \times \mathbb{M}^{3 \times 3} \rightarrow \mathbb{R}$ is a Carath\'eodory function satisfying conditions of $p$-coercivity and  $p$-growth : for some $0< \beta' \leq \beta <+\infty$ and some  $1<p<\infty$,
\begin{equation}\label{A}
\beta ' |F|^p \leq W(x;F) \leq \beta(|F|^p+1),\; F \in \mathbb{M}^{3 \times 3},\; \mbox{ for a.e. } x \in \Omega.
\end{equation}
The usual Euclidian norm on the space $\mathbb{M}^{N \times m}$ of real $N \times m$ matrices is denoted by $|F|$. The minimisation problem becomes
\begin{equation}\label{C}
\inf_{v-x_{\varepsilon} \in W^{1,p}_{\Gamma}(\Omega;\mathbb{R}^3)} \mathcal{E}_{\varepsilon}(v),
\end{equation}
where $W^{1,p}_{\Gamma}(\Omega;\mathbb{R}^3)$ stands for the functions in $W^{1,p}(\Omega;\mathbb{R}^3)$ with zero trace on the lateral boundary  $\Gamma:=\partial \omega \times (-1,1)$.

If we denote by $ g_0 ^\pm$ (resp. $g^\pm$) the trace of $ g_0$ (resp. $g$) on $\omega \times \{\pm 1\}$, in view of Remark 3.2.3 of \cite{Fox&Ra&Si}, the loading vectors $ g_0^+$ and $ g_0^-$ must satisfy $ g_0^+ + g_0^- =0$. In the second section, we assume the stronger condition that $ g_0^+= g_0^-=0$ . The physical implication of this assumption is that the plate of thickness $2\varepsilon$ cannot support a non vanishing resultant surface load as the thickness $\varepsilon$ goes to zero. We generalize here the result of \cite{Le&Ra} and \cite{Bra&Fo&Fr} to a general inhomogeneity. In the third section, we address the general case of admissible surface loadings. It deals with a similar problem, in which the class of surface forces generates a bending moment density as in \cite{Bo&Fo&Ma}; the limit behavior is not solely characterized by the limit deformations (a $\mathbb{R}^3$-valued field defined on the mid-plane), but it also involves the average of the Cosserat vector also defined on the mid-plane. Once again, we generalize  the result of \cite{Bo&Fo&Ma} to the inhomogeneous case. The fourth and last section demonstrates that the classical membrane model can be seen as a particular case of the Cosserat model when the bending moment density is zero.

As for notation,  $\mathcal{A}(\omega)$ is the family of open subsets of $\omega$; $\mathcal{L}^N$ stands for the N-dimensional Lebesgue measure in $\mathbb{R}^N$ (N=2 or 3); $\rightarrow$  always denotes strong convergence whereas $\rightharpoonup$ (resp. $\stackrel{*}{\rightharpoonup}$) denotes weak (resp. weak-*) convergence. Finally, we loosely identify $L^p(\omega;\mathbb{R}^3)$ (resp. $W^{1,p}(\omega;\mathbb{R}^3)$) with those functions in $L^p(\Omega;\mathbb{R}^3)$ (resp. $W^{1,p}(\Omega;\mathbb{R}^3)$) that do not depend upon $x_3$.

\section{Classical nonlinear membrane model}

In this section, we assume that $g_{\varepsilon}= \varepsilon g$ with $g \in L^{p'}(\Sigma;\mathbb{R}^3)$. Thus, the minimisation problem (\ref{C}) becomes
$$\inf_{v-x_{\varepsilon} \in W^{1,p}_{\Gamma}(\Omega;\mathbb{R}^3)} \left\{ \int_{\Omega}W\left(x;D_{\alpha} v \Big| \frac{1}{\varepsilon} D_3 v\right) dx - \int_{\Omega}f.v \, dx -\int_{\Sigma}g.v \, d\mathcal H^2 \right\}.$$
Define for any $(u;A) \in L^p(\Omega;\mathbb{R}^3) \times \mathcal{A}(\omega)$,
$$J_{\varepsilon}(u;A):=\left\{
\begin{array}{ll}
\ds\int_{A \times (-1,1)}W\left(x;D_{\alpha} u\Big|\frac{1}{\varepsilon}D_3 u\right) dx & \text{ if } u \in W^{1,p}(A \times (-1,1);\mathbb{R}^3),\\\\
+\infty & \text{ otherwise},
\end{array}
\right.$$
and
\begin{equation}\label{F}
J_{\{\varepsilon\}}(u;A):=\inf_{\{u_{\varepsilon}\}}\left\{ \liminf_{\varepsilon \rightarrow 0} J_{\varepsilon}(u_{\varepsilon};A) :  u_{\varepsilon} \rightarrow u \text{ in }  L^p(A \times (-1,1);\mathbb{R}^3) \right\}.
\end{equation}
\begin{rem}\label{dem2}
 For any $A \in \mathcal{A}(\omega)$,  $J_{\{\varepsilon\}}(u;A)=+\infty$
 whenever $u \in L^p(\Omega;\mathbb{R}^3) \setminus W^{1,p}(A;\mathbb{R}^3)$, as is easily seen in view of the definition of
 $J_{\varepsilon}$, together with the coercivity condition (\ref{A}) .
\end{rem}
By virtue of Remark \ref{dem2}, together with Theorem 2.5 in \cite{Bra&Fo&Fr}, for all sequences $\{\varepsilon\}$, there exists a subsequence $\{\varepsilon_n\}$ such that $J_{\{\varepsilon_n\}}(.;A)$ defined in (\ref{F}) is the $\Gamma(L^p)-$limit of $J_{\varepsilon_n}(.;A)$. Further, there exists a Carath\'eodory function $W_{\{\varepsilon_n\}}:\omega \times \mathbb{M}^{3 \times 2} \rightarrow \mathbb{R}$ such that, for all $A \in \mathcal{A}(\omega)$ and all $u \in W^{1,p}(A;\mathbb{R}^3)$
$$J_{\{\varepsilon_n\}}(u;A)=2\int_A W_{\{\varepsilon_n\}}(x_{\alpha};D_{\alpha} u)dx_{\alpha}.$$
\begin{rem}
\label{bci} Lemma 2.6 of \cite{Bra&Fo&Fr} implies that $J_{\{\varepsilon_n\}}(u;A)$ is unchanged if the approximating sequences
$\{u_{\varepsilon_n}\}$ are constrained to match the lateral boundary condition of their target, i.e. $u_{\varepsilon_n}\equiv u \mbox{ on }\partial A \times (-1,1)$.
\end{rem}

From now onward, we will assume that $\{\varepsilon_n\}$ denotes a subsequence of $\{\varepsilon\}$ such that the $\Gamma(L^p)$-limit of $J_{\varepsilon_n}(u;A)$ exists,  in which case it coincides with $J_{\{\varepsilon_n\}}(u;A)$.  Under the hypothesis that $W$ is a homogeneous elastic energy density, it is proved in \cite{Le&Ra}, Theorem 2, that  $J_{\{\varepsilon_n\}}(u;A)$ does not depend upon the choice of the sequence $\{\varepsilon_n\}$. It is given by
$$J_{\{\varepsilon_n\}}(u;A)=2\int_A \mathcal{Q}\overline{W}(D_{\alpha} u)dx_{\alpha},$$
with for all $\overline{F} \in \mathbb{M}^{3 \times 2}$,
$$\overline{W}(\overline{F}):= \inf_{z \in \mathbb{R}^3}W(\overline{F}|z),$$
and,
$$\mathcal{Q}\overline{W}(\overline{F}):= \inf_{\phi \in W^{1,p}_0(Q';\mathbb{R}^3)} \int_{Q'}\overline{W}(\overline{F} + D_{\alpha}\phi)dx_{\alpha},$$
where $Q':=(0,1)^2$, and $\mathcal{Q}\overline{W}$ is the  2D-quasiconvexification of $\overline{W}$. This result was extended to the case where $W$ is also function of $x_3$ in \cite{Bra&Fo&Fr}, Theorem \ref{t3.1}. It is proved there that, in such a case, $J_{\{\varepsilon_n\}}$ is given by
$$J_{\{\varepsilon_n\}}(u;A)=2\int_A \underline{W}(D_{\alpha} u)dx_{\alpha},$$
with for all $\overline{F} \in \mathbb{M}^{3 \times 2}$,
\begin{eqnarray*}
\underline{W}(\overline{F}) & := & \inf_{L,\phi} \Big\{ \frac{1}{2}\int_{Q' \times (-1,1)} W(x_3;\overline{F}+D_{\alpha} \phi|L D_3 \phi)dx_{\alpha}dx_3 : \; L>0,\\
 & & \phi \in W^{1,p}(Q' \times (-1,1);\mathbb{R}^3), \, \phi=0 \text{ on }\partial Q' \times (-1,1) \Big\}.
\end{eqnarray*}
We wish to extend those results to the case where $W$ is a function of {\bf both} $x_3$ {\bf and} $x_{\alpha}$. We set,  for all $\overline{F} \in \mathbb{M}^{3 \times 2}$ and for a.e. $x_0 \in \omega$,
\begin{eqnarray}
\underline{W}(x_0;\overline{F}) & := & \inf_{L,\phi} \Big\{ \frac{1}{2}\int_{Q' \times (-1,1)} W(x_0,x_3;\overline{F}+D_{\alpha} \phi|L D_3 \phi)dx_{\alpha}dx_3 : \; L>0,\nonumber\\
 & & \phi \in W^{1,p}(Q' \times (-1,1);\mathbb{R}^3), \, \phi=0 \text{ on }\partial Q' \times (-1,1) \Big\}.\label{D}
\end{eqnarray}
The following theorem holds :
\begin{theo}
\label{t2.1}
For all $\overline{F} \in \mathbb{M}^{3 \times 2}$ and for a.e. $x_0 \in \omega$, $W_{\{\varepsilon_n\}}(x_0;\overline{F})=\underline{W}(x_0;\overline{F})$ where $\underline{W}$ is given by (\ref{D}). Furthermore, for any  $ A\in \mathcal{A}(\omega)$, $J_{\varepsilon}(.;A)$ $\Gamma(L^p)$-converges to $J_{\{\varepsilon\}}(.;A)$ and for all $u \in L^p(\Omega;\mathbb{R}^3)$,
$$J_{\{\varepsilon\}}(u;A)=\left\{
\begin{array}{ll}
\ds2\int_A \underline{W}(x_{\alpha};D_{\alpha} u)\,dx_{\alpha} & \text{ if } u \in W^{1,p}(A;\mathbb{R}^3),\\\\
+\infty & \text{otherwise}.
\end{array}
\right.$$
\end{theo}
The proof of this Theorem is a direct consequence of  Lemmata \ref{lem>},
\ref{lem<} below.
\begin{lem}\label{lem>}
For all $\overline{F} \in \mathbb{M}^{3 \times 2}$ and for a.e. $x_0 \in \omega$,
$$W_{\{\varepsilon_n\}}(x_0,\overline{F}) \geq \underline{W}(x_0,\overline{F}).$$
\end{lem}
\textit{Proof. }Let us fix $\overline F \in \mathbb M^{3 \times 2}$, we set $u(x_\alpha):= \overline F x_\alpha$ and let $x_0$ be a Lebesgue point of both $W_{\{\varepsilon_n\}}(.;\overline F)$ and $\underline W(.;\overline F)$. We denote by $Q'(x_0,r)$, the cube of $\mathbb R^2$ of center $x_0$ and side length $r$, where $r>0$ is fixed and small enough so that $Q'(x_0,r) \in \mathcal A(\omega)$. According to the equi-integrability Theorem (Theorem 1.1 in \cite{Bo&Fo}), there exists a subsequence of $\{\varepsilon_n\}$ (not relabelled) and a sequence $\{u_n\} \subset W^{1,p}(Q'(x_0,r) \times (-1,1);\mathbb R^3)$ such that
$$\left\{ \begin{array}{l}
u_n \to 0 \text{ in }L^p(Q'(x_0,r) \times (-1,1);\mathbb R^3),\\\\
\left\{ \left| \left( D_\alpha u_n \big| \frac{1}{\varepsilon_n} D_3 u_n\right) \right|^p \right\} \text{ is equi-integrable},\\\\
\ds J_{\{\varepsilon_n\}}(u;Q'(x_0,r)) = \lim_{n \to +\infty} \int_{Q'(x_0,r) \times (-1,1)}W\left(x_\alpha,x_3;\overline F +D_\alpha u_n \Big| \frac{1}{\varepsilon_n} D_3 u_n\right)dx_\alpha dx_3.
\end{array} \right.$$
Set $$F_n(x):=\left(\overline F +D_\alpha u_n(x) \Big| \frac{1}{\varepsilon_n} D_3 u_n(x)\right).$$
For any $h \in \mathbb N$, we cover $Q'(x_0,r)$ with $h^2$ disjoints cubes $Q'_{i,h}$ of side length $r/h$. Thus $Q'(x_0,r)=\bigcup_{i=1}^{h^2} Q'_{i,h}$ and
\begin{equation}\label{equiint}
J_{\{\varepsilon_n\}}(u;Q'(x_0,r)) = \left( \limsup_{h \to +\infty} \right) \limsup_{n \to +\infty} \sum_{i=1}^{h^2}\int_{Q'_{i,h} \times (-1,1)}W(x;F_n(x))dx.
\end{equation}
Since $W$ is a Carath\'eodory integrand, Scorza-Dragoni's Theorem (see \cite{Ek&Te}, Chapter VIII) implies the existence, for any $\eta >0$, of a compact set $K_\eta \subset \Omega$ such that
\begin{equation}\label{mes}
\mathcal L^3(\Omega \setminus K_\eta) < \eta,
\end{equation}
and the restriction of $W$ to $K_\eta \times \mathbb M^{3 \times 3}$ is continuous. For any $\lambda>0$, define
$$R^\lambda_n:=\left\{ x \in Q'(x_0,r) \times (-1,1): |F_n(x)| \leq \lambda\right\}.$$
By virtue of Chebyshev's inequality, there exists a constant $C>0$ -- which does not depend on $n$ or $\lambda$ -- such that
\begin{equation}\label{lambda}
\mathcal L^3 ([Q'(x_0,r) \times (-1,1)] \setminus R^\lambda_n) < \frac{C}{\lambda^p}.
\end{equation}
Denoting by $W^{\eta,\lambda}$ the continuous extension  of $W$ outside $K_\eta \times \overline B(0,\lambda)$ (defined e.g. in Theorem 1, Section 1.2 in \cite{Ev&Ga}), $W^{\eta,\lambda}$ is continuous on $\mathbb R^3 \times \mathbb M^{3 \times 3}$ and satisfies the following bound
\begin{equation}\label{ext}
0 \leq W^{\eta,\lambda}(x;F) \leq \max_{K_\eta \times \overline B(0,\lambda)}W \leq \beta(1+\lambda^p) \quad \text{for all } (x;F) \in \mathbb R^3 \times \mathbb M^{3 \times 3}.
\end{equation}
In view of (\ref{equiint}), we have
$$J_{\{\varepsilon_n\}}(u;Q'(x_0,r)) \geq \limsup_{\lambda \to +\infty} \limsup_{\eta \to 0} \limsup_{h \to +\infty} \limsup_{n \to +\infty} \sum_{i=1}^{h^2} \int_{[Q'_{i,h} \times (-1,1)] \cap R^\lambda_n \cap K_\eta}W^{\eta,\lambda} (x;F_n(x))dx.$$
By virtue of (\ref{ext}) and (\ref{mes}),
$$\sum_{i=1}^{h^2} \int_{[Q'_{i,h} \times (-1,1)] \cap R^\lambda_n \setminus K_\eta}W^{\eta,\lambda} (x;F_n(x))dx \leq \beta(1+\lambda^p)\eta \xrightarrow[\eta \to 0]{}0,$$
uniformly in $(n,h)$. Therefore
$$J_{\{\varepsilon_n\}}(u;Q'(x_0,r)) \geq \limsup_{\lambda \to +\infty} \limsup_{\eta \to 0} \limsup_{h \to +\infty} \limsup_{n \to +\infty} \sum_{i=1}^{h^2} \int_{[Q'_{i,h} \times (-1,1)] \cap R^\lambda_n}W^{\eta,\lambda} (x;F_n(x))dx.$$
Since $W^{\eta,\lambda}$ is continuous, it is uniformly continuous on $\overline \Omega \times \overline B(0,\lambda)$. Thus there exists a continuous and increasing function $\omega_{\eta,\lambda} : \mathbb R^+ \longrightarrow \mathbb R^+$ satisfying $\omega_{\eta,\lambda}(0)=0$ and such that
\begin{equation}\label{cont}
|W^{\eta,\lambda}(x;F)-W^{\eta,\lambda}(y;G)| \leq \omega_{\eta,\lambda}(|x-y|+|F-G|), \quad \forall \, (x;F),\; (y;G) \in \overline \Omega \times \overline B(0,\lambda).
\end{equation}
Consenquently, for all $(x_\alpha,x_3) \in  [Q'_{i,h} \times (-1,1)] \cap R^\lambda_n $ and all $y_\alpha \in Q'_{i,h}$,
$$\big|W^{\eta,\lambda}(x_\alpha,x_3;F_n(x_\alpha,x_3))-W^{\eta,\lambda}(y_\alpha,x_3;F_n(x_\alpha,x_3))\big| \leq  \omega_{\eta,\lambda}(|x_\alpha-y_\alpha|) \leq \omega_{\eta,\lambda}(c/h).$$
We get, after integration in $(x,y_\alpha)$ and summation,
\begin{eqnarray*}
 & \ds  \sup_{n \in \mathbb N}\sum_{i=1}^{h^2} \frac{h^2}{r^2} \int_{Q'_{i,h}} \left\{ \int_{R^\lambda_n \cap [ Q'_{i,h} \times (-1,1)]}\big|W^{\eta,\lambda}(y_\alpha,x_3;F_n(x))-W^{\eta,\lambda}(x_\alpha,x_3;F_n(x))\big|dx \right\} dy_\alpha\\
 & \ds \leq 2 r^2\omega_{\eta,\lambda}(c/h) \xrightarrow[h \to +\infty]{} 0.
\end{eqnarray*}
Hence,
\begin{eqnarray*}
&\ds J_{\{\varepsilon_n\}}(u;Q'(x_0,r)) \geq \\
&\ds \limsup_{\lambda \to +\infty} \limsup_{\eta \to 0} \limsup_{h \to +\infty} \limsup_{n \to +\infty} \sum_{i=1}^{h^2} \frac{h^2}{r^2} \int_{ Q'_{i,h}} \left\{ \int_{[ Q'_{i,h} \times (-1,1)] \cap R^\lambda_n }W^{\eta,\lambda}(y_\alpha,x_3;F_n(x))dx \right\} dy_\alpha.
\end{eqnarray*}
Define the following sets which depend on all parameters $(\eta,\lambda,i,h,n)$ :
\begin{eqnarray*}
&&E := \{ (y_\alpha,x_\alpha,x_3) \in Q'_{i,h}\times Q'_{i,h} \times (-1,1) : (y_\alpha,x_3) \in K_\eta \text{ and } (x_\alpha,x_3) \in R^\lambda_n\},\\
&&E_1 := \{ (y_\alpha,x_\alpha,x_3) \in Q'_{i,h}\times Q'_{i,h} \times (-1,1) : (y_\alpha,x_3) \not \in K_\eta \text{ and } (x_\alpha,x_3) \in R^\lambda_n\},\\
&&E_2:= \{ (y_\alpha,x_\alpha,x_3) \in Q'_{i,h}\times Q'_{i,h} \times (-1,1) : (x_\alpha,x_3) \not \in R^\lambda_n\},
\end{eqnarray*}
and note that $Q'_{i,h} \times Q'_{i,h} \times (-1,1)= E \cup E_1 \cup E_2$. Since $W$ and $W^{\eta,\lambda}$ coincide on $K_\eta \times \overline B(0,\lambda)$,
\begin{eqnarray}
J_{\{\varepsilon_n\}}(u;Q'(x_0,r)) & \geq  & \limsup_{\lambda \to +\infty} \limsup_{\eta \to 0} \limsup_{h \to +\infty} \limsup_{n \to +\infty} \sum_{i=1}^{h^2} \frac{h^2}{r^2} \int_{E} W^{\eta,\lambda}(y_\alpha,x_3;F_n(x))dx \, dy_\alpha \nonumber\\
 & = & \limsup_{\lambda \to +\infty} \limsup_{\eta \to 0} \limsup_{h \to +\infty} \limsup_{n \to +\infty} \sum_{i=1}^{h^2} \frac{h^2}{r^2} \int_{E} W(y_\alpha,x_3;F_n(x))dx \, dy_\alpha.\label{E}
\end{eqnarray}
We will prove that the corresponding terms over $E_1$ and $E_2$ are zero. Indeed, in view of (\ref{mes}) and the $p$-growth condition (\ref{A}),
\begin{eqnarray}
\sum_{i=1}^{h^2}  \frac{h^2}{r^2}\int_{E_1} W(y_\alpha,x_3;F_n(x))dx  dy_\alpha & \leq & \sum_{i=1}^{h^2}  \frac{h^2}{r^2} \, \mathcal L^2(Q'_{i,h}) \, \mathcal L^3([Q'_{i,h} \times (-1,1)] \setminus K_\eta)\, \beta (1+\lambda^p)\nonumber\\
 & = & \beta (1+\lambda^p)\, \mathcal L^3([Q'(x_0,r) \times (-1,1)] \setminus K_\eta)\nonumber\\
 & < & \beta (1+\lambda^p)\eta \xrightarrow[\eta \to 0]{} 0,\label{E1}
\end{eqnarray}
uniformly in $(n,h)$. The bound from above in (\ref{A}), the equi-integrability of $\{|F_n|^p\}$ and (\ref{lambda}) imply that
\begin{eqnarray}
\sum_{i=1}^{h^2}  \frac{h^2}{r^2}\int_{E_2} W(y_\alpha,x_3;F_n(x))dx  dy_\alpha & \leq & \sum_{i=1}^{h^2}  \frac{h^2}{r^2} \, \mathcal L^2(Q'_{i,h}) \, \beta \int_{[Q'_{i,h} \times (-1,1)] \setminus R^\lambda_n}(1+|F_n(x)|^p)dx\nonumber\\
 & = & \beta \int_{[Q'(x_0,r) \times (-1,1)] \setminus R^\lambda_n}(1+|F_n(x)|^p)dx \xrightarrow[\lambda \to +\infty]{} 0,\label{E2}
\end{eqnarray}
uniformly in $(\eta,n,h)$. Thus, in view of (\ref{E}), (\ref{E1}), (\ref{E2}), Fatou's Lemma yields
\begin{eqnarray*}
J_{\{\varepsilon_n\}}(u;Q'(x_0,r)) & \geq  & \limsup_{h \to +\infty} \limsup_{n \to +\infty} \sum_{i=1}^{h^2}  \frac{h^2}{r^2}\int_{Q'_{i,h}} \left\{ \int_{Q'_{i,h} \times (-1,1)}W(y_\alpha,x_3;F_n(x))dx \right\} dy_\alpha\\
 & \geq  & \limsup_{h \to +\infty} \liminf_{n \to +\infty} \sum_{i=1}^{h^2}  \frac{h^2}{r^2}\int_{Q'_{i,h}} \left\{ \int_{Q'_{i,h} \times (-1,1)}W(y_\alpha,x_3;F_n(x))dx \right\} dy_\alpha\\
 & \geq & \limsup_{h \to +\infty} \sum_{i=1}^{h^2}  \frac{h^2}{r^2} \int_{Q'_{i,h}} \left\{ \liminf_{n \to +\infty}\int_{Q'_{i,h} \times (-1,1)}W(y_\alpha,x_3;F_n(x))dx \right\} dy_\alpha,
\end{eqnarray*}
We apply, for a.e. $y_\alpha \in Q'_{i,h}$, Theorem 3.1 in \cite{Bra&Fo&Fr} to the Carath\'eodory function $(x_3;F) \mapsto W(y_\alpha,x_3;F)$; in particular
$$\liminf_{n \to +\infty}\int_{Q'_{i,h} \times (-1,1)}W(y_\alpha,x_3;F_n(x))dx \geq \frac{2 r^2}{h^2}\; \underline W(y_\alpha;\overline F).$$
Thus
$$J_{\{\varepsilon_n\}}(u;Q'(x_0,r))  \geq   \limsup_{h \to +\infty}\sum_{i=1}^{h^2} \frac{h^2}{r^2} \int_{Q'_{i,h}} \frac{2 r^2}{h^2}\; \underline W(y_\alpha;\overline F) dy_\alpha = 2 \int_{Q'(x_0,r)} \underline W(y_\alpha;\overline F) dy_\alpha.$$
Dividing both sides of the previous inequality by $r^2$ and passing to the limit when $r \searrow 0^+$, we obtain
$$W_{\{\varepsilon_n\}}(x_0;\overline F) \geq  \underline W (x_0;\overline F).$$
\hfill $\Box$
\begin{lem}\label{lem<}
For all $\overline{F} \in \mathbb{M}^{3 \times 2}$ and for a.e. $x_0 \in \omega$,
$$W_{\{\varepsilon_{n}\}}(x_0;\overline{F}) \leq \underline{W}(x_0;\overline{F}).$$
\end{lem}
\textit{Proof. }For all $k \geq 1$, let $L_k>0$ and $\varphi_k \in W^{1,\infty}(Q' \times (-1,1);\mathbb R^3)$ with $\varphi_k=0$ on $\partial Q' \times (-1,1)$ be such that
\begin{equation}\label{fk}
Z_k(x_0;\overline F):=\frac{1}{2}\int_{Q' \times (-1,1)}W(x_0,x_3; \overline F +D_\alpha \varphi_k|L_k D_3 \varphi_k)dx_\alpha dx_3 \leq \underline W(x_0;\overline F)+\frac{1}{k}.
\end{equation}
This is legitimate because of the density of $W^{1,\infty}(Q' \times (-1,1);\mathbb R^3)$ into $W^{1,p}(Q' \times (-1,1);\mathbb R^3)$ and the $p$-growth condition (\ref{A}). We extend $\varphi_k$ to $\mathbb R^2 \times (-1,1)$ by $Q'$-periodicity and set $F_k(x):= (\overline F +D_\alpha \varphi_k(x)|L_k D_3 \varphi_k(x))$. Then, there exists $M_k >0$ such that
\begin{equation}\label{borne}
\| F_k\|_{L^\infty(\mathbb R^2 \times (-1,1);\mathbb R^3)}\leq M_k.
\end{equation}
Let $\overline F \in \mathbb M^{3 \times 2}$ and $x_0$ be a Lebesgue point of $\underline W (.;\overline F)$ and $Z_k(.;\overline F)$ for all $k\geq 1$. We choose $r>0$ small enough such that $Q'(x_0,r) \in \mathcal A(\omega)$. Fix $k \geq 1$ and set
$$\left\{
\begin{array}{l}
\ds u(x_\alpha)  :=  \overline F x_\alpha,\\
\ds u_n^k(x_\alpha,x_3)  :=  \overline F x_\alpha+L_k \varepsilon_n \varphi_k\left( \frac{x_\alpha}{L_k \varepsilon_n},x_3 \right).
\end{array}
\right.$$
Since $u_n^k \xrightarrow[n \to +\infty]{} u$ in $L^p(Q'(x_0,r) \times (-1,1);\mathbb R^3)$,
\begin{eqnarray*}
J_{\{\varepsilon_n\}}(u;Q'(x_0,r)) & \leq & \liminf_{n \to +\infty} \int_{Q'(x_0,r) \times (-1,1)}W\left(x_\alpha,x_3;D_\alpha u_n^k \Big| \frac{1}{\varepsilon_n}D_3 u_n^k \right)dx_\alpha x_3\\
 & = & \liminf_{n \to +\infty} \int_{Q'(x_0,r) \times (-1,1)}W\left(x_\alpha,x_3;F_k\left(\frac{x_\alpha}{L_k \varepsilon_n},x_3 \right) \right)dx_\alpha x_3.
\end{eqnarray*}
As before, we split $Q'(x_0,r)$ into $h^2$ disjoint cubes $Q'_{i,h}$ of length $r/h$. Then,
$$J_{\{\varepsilon_n\}}(u;Q'(x_0,r)) \leq \left( \liminf_{h \to +\infty} \right) \liminf_{n \to +\infty} \sum_{i=1}^{h^2}\int_{Q'_{i,h} \times (-1,1)}W\left(x_\alpha,x_3;F_k\left(\frac{x_\alpha}{L_k \varepsilon_n},x_3 \right) \right)dx_\alpha x_3.$$
Let $K_\eta$ be like in Lemma \ref{lem>} and $W^{\eta,k}$ be a continuous extension of $W$ outside $K_\eta \times \overline B(0,M_k)$ which satisfies the analogue of (\ref{ext}) with $M_k$ instead of $\lambda$. In view of the $p$-growth condition (\ref{A}), (\ref{borne}) and (\ref{mes}), we get
$$\sup_{n \in \mathbb N} \sum_{i=1}^{h^2} \int_{[Q'_{i,h}  \times (-1,1)] \setminus K_\eta}W\left(x_\alpha,x_3;F_k\left(\frac{x_\alpha}{L_k \varepsilon_n},x_3 \right) \right)dx_\alpha x_3 \leq \beta(1+M_k^p)\eta \xrightarrow[\eta \to 0]{} 0.$$
Thus,
\begin{eqnarray*}
J_{\{\varepsilon_n\}}(u;Q'(x_0,r)) & \leq & \liminf_{\eta \to 0} \liminf_{h \to +\infty} \liminf_{n \to +\infty} \sum_{i=1}^{h^2}\int_{[Q'_{i,h}  \times (-1,1)] \cap K_\eta}W^{\eta,k} \left(x_\alpha,x_3;F_k\left(\frac{x_\alpha}{L_k \varepsilon_n},x_3 \right) \right)dx_\alpha x_3\\
& \leq & \liminf_{\eta \to 0} \liminf_{h \to +\infty} \liminf_{n \to +\infty} \sum_{i=1}^{h^2}\int_{Q'_{i,h}  \times (-1,1)}W^{\eta,k} \left(x_\alpha,x_3;F_k\left(\frac{x_\alpha}{L_k \varepsilon_n},x_3 \right) \right)dx_\alpha x_3.
\end{eqnarray*}
Since $W^{\eta,k}$ is continuous, it is uniformly continuous on $\overline \Omega \times \overline B(0,M_k)$. Thus, there exists a continuous and increasing function $\omega_{\eta,k} : \mathbb R^+ \longrightarrow \mathbb R^+$ satisfying $\omega_{\eta,k}(0)=0$ and the analogue of (\ref{cont}), replacing  $\lambda$ by $M_k$. Then, for every $(x_\alpha,x_3) \in Q'_{i,h} \times (-1,1)$ and every $y_\alpha \in Q'_{i,h}$,
\begin{eqnarray*}
\left|W^{\eta,k} \left(x_\alpha,x_3;F_k\left(\frac{x_\alpha}{L_k \varepsilon_n},x_3 \right) \right)-W^{\eta,k} \left(y_\alpha,x_3;F_k\left(\frac{x_\alpha}{L_k \varepsilon_n},x_3 \right) \right)\right| & \leq & \omega_{\eta,k}(|x_\alpha-y_\alpha|)\\
 & \leq &  \omega_{\eta,k}(c/h).
\end{eqnarray*}
Integration and summation yield in turn
\begin{eqnarray*}
 & \ds \sup_{n \in \mathbb N}\sum_{i=1}^{h^2} \frac{h^2}{r^2} \int_{Q'_{i,h}} \left\{ \int_{Q'_{i,h} \times (-1,1)}\left|W^{\eta,k}(y_\alpha,x_3;F_k\left(\frac{x_\alpha}{L_k \varepsilon_n},x_3 \right)-W^{\eta,k}(x_\alpha,x_3;F_k\left(\frac{x_\alpha}{L_k \varepsilon_n},x_3 \right)\right|dx \right\} dy_\alpha\\
& \ds \leq 2 r^2 \omega_{\eta,k}(c/h) \xrightarrow[h \to +\infty]{} 0.
\end{eqnarray*}
Hence,
\begin{eqnarray*}
& J_{\{\varepsilon_n\}}(u;Q'(x_0,r)) \leq\\
& \ds \liminf_{\eta \to 0} \liminf_{h \to +\infty} \liminf_{n \to +\infty} \sum_{i=1}^{h^2} \frac{h^2}{r^2} \int_{Q'_{i,h}} \left\{ \int_{Q'_{i,h} \times (-1,1)} W^{\eta,k}(y_\alpha,x_3;F_k\left(\frac{x_\alpha}{L_k \varepsilon_n},x_3 \right)dx_\alpha dx_3 \right\} dy_\alpha.
\end{eqnarray*}
According to (\ref{ext}) and (\ref{mes}),
$$\sup_{n \in \mathbb N}\sum_{i=1}^{h^2} \frac{h^2}{r^2} \int_{Q'_{i,h}} \left\{ \int_{[Q'_{i,h} \times (-1,1)] \setminus K_\eta} W^{\eta,k}(y_\alpha,x_3;F_k\left(\frac{x_\alpha}{L_k \varepsilon_n},x_3 \right)dy_\alpha dx_3 \right\} dx_\alpha \leq \beta(1+M_k^p)\eta \xrightarrow[\eta \to 0]{} 0.$$
Since $W^{\eta,k}$ coincides with $W$ on $K_\eta \times \overline B(0,M_k)$, we get
\begin{eqnarray*}
& J_{\{\varepsilon_n\}}(u;Q'(x_0,r))\\
&\leq \ds \liminf_{\eta \to 0} \liminf_{h \to +\infty} \liminf_{n \to +\infty} \sum_{i=1}^{h^2} \frac{h^2}{r^2} \int_{Q'_{i,h}} \left\{ \int_{[Q'_{i,h} \times (-1,1)] \cap K_\eta} W(y_\alpha,x_3;F_k\left(\frac{x_\alpha}{L_k \varepsilon_n},x_3 \right)dy_\alpha dx_3 \right\} dx_\alpha\\
&\leq \ds  \liminf_{h \to +\infty} \sum_{i=1}^{h^2} \frac{h^2}{r^2}\limsup_{n \to +\infty} \int_{Q'_{i,h}} \left\{ \int_{Q'_{i,h} \times (-1,1)} W(y_\alpha,x_3;F_k\left(\frac{x_\alpha}{L_k \varepsilon_n},x_3 \right)dy_\alpha dx_3 \right\} dx_\alpha.
\end{eqnarray*}
Riemann-Lebesgue's Lemma applied to the $Q'$-periodic function $\int_{Q'_{i,h} \times (-1,1)} W (y_\alpha,x_3;F_k(.,x_3) )dy_\alpha dx_3$ implies that
$$J_{\{\varepsilon_n\}}(u;Q'(x_0,r)) \leq  \liminf_{h \to +\infty}\sum_{i=1}^{h^2} \frac{h^2}{r^2}\int_{Q'_{i,h}} \frac{2 r^2}{h^2}Z_k(y_\alpha;\overline F)dy_\alpha = 2 \int_{Q'(x_0,r)}Z_k(y_\alpha;\overline F)dy_\alpha.$$
Dividing both sides of the inequality by $r^2$ and letting $r \searrow 0^+$, we get in view of the definition of $x_0$ and (\ref{fk}),
$$W_{\{\varepsilon_n\}}(x_0;\overline F) \leq Z_k(x_0;\overline F) \leq \underline W (x_0;\overline F)+\frac{1}{k}.$$
Passing to the limit when $k \nearrow +\infty$ yields the desired result.
\hfill$\Box$
\\
\\
\textit{Proof of Theorem \ref{t2.1}.} For a.e. $x_0 \in \omega$ and for all $\overline{F} \in \mathbb{M}^{3 \times 2}$, $W_{\{\varepsilon_n\}}(x_0;\overline{F}) = \underline{W}(x_0;\overline{F})$. Since  the $\Gamma(L^p)$-limit does not depend upon the choice of sequence $\{\varepsilon_n\}$, appealing to Proposition 7.11  in \cite{Bra&De}  we conclude that for any $A \in \mathcal{A}(\omega)$, the whole sequence $J_{\varepsilon}(.;A)$ $\Gamma(L^p)$-converges to $J_{\{\varepsilon\}}(.;A)$ and we have,
$$J_{\{\varepsilon\}}(u;A)=2\int_A \underline{W}(x_{\alpha};D_{\alpha}u)dx_{\alpha},$$
 for all $u \in W^{1,p}(A;\mathbb{R}^3)$.\hfill$\Box$
\begin{rem}
 Proposition 4.1 gives another expression for the energy density $\underline{W}$.
\end{rem}
\begin{rem}
By construction and thanks to Remark 3.3 of \cite{Bra&Fo&Fr}, Theorem \ref{t2.1} generalizes both Theorem 2 of \cite{Le&Ra} and Theorem \ref{t3.1} of \cite{Bra&Fo&Fr}.
\end{rem}

\section{Cosserat nonlinear membrane model}

In this section, we assume as in \cite{Bo&Fo&Ma} that $g_{\varepsilon}:= g_0+\varepsilon g$ with $g_0, g \in L^{p'}(\Sigma;\mathbb{R}^3)$ and $g_0^+ + g_0^- =0$. Thus, the minimisation problem (\ref{C}) reads as
$$\inf_{v-x_{\varepsilon} \in W^{1,p}_{\Gamma}(\Omega;\mathbb{R}^3)} \left\{ \int_{\Omega}W\left(x;D_{\alpha} v \Big| \frac{1}{\varepsilon} D_3 v\right) dx - L_{\varepsilon}(v) \right\},$$
with
$$L_{\varepsilon}(v):=\int_{\Omega}f.v \, dx +\int_{\Sigma}g.v \, d\mathcal H^2 + \int_{\omega}g_0^+ . \left( \frac{v^+-v^-}{\varepsilon}\right)dx_{\alpha}, \quad v^\pm(x_\alpha):=v(x_\alpha,\pm 1).$$
If $v_{\varepsilon} \to v$ in $L^p(\Omega;\mathbb{R}^3)$ is a minimizing sequence and if $ b_{\varepsilon}:=\frac{1}{\varepsilon}D_3 v_{\varepsilon}$, then
$$L_{\varepsilon}(v_{\varepsilon})=\int_{\Omega}f.v_{\varepsilon} \, dx +\int_{\Sigma} g.v_\varepsilon \, d\mathcal H^2+ 2\int_{\omega}g_0^+.\overline{b}_{\varepsilon}dx_{\alpha}, \quad \text{where } \overline b_\varepsilon = \frac{1}{2} \int_{-1}^1 b_\varepsilon(.,x_3) dx_3.$$
By virtue of the coercivity condition (\ref{A}), we deduce that the sequence $\{v_{\varepsilon}\}$ is uniformly bounded in $W^{1,p}(\Omega;\mathbb{R}^3)$ and that, for a subsequence of $\{\varepsilon\}$ still labelled $\{\varepsilon\}$, $v_{\varepsilon} \rightharpoonup v$ in $W^{1,p}(\Omega;\mathbb{R}^3)$ and $b_{\varepsilon} \rightharpoonup b$ in $L^p(\Omega;\mathbb{R}^3)$ with $v \in W^{1,p}(\omega;\mathbb{R}^3)$. As previously, $v$ is associated to the mid-plane deformation,  whereas $b$ is the Cosserat vector. In any case, $\lim_{\varepsilon \rightarrow 0} L_{\varepsilon}(v_{\varepsilon}) = L(v,\overline{b})$, with
\begin{equation}\label{bm:0}
L(v,\overline{b}):=\int_{\omega}\left( 2 \overline{f} + g^+ +g^- \right).v \, dx_{\alpha} + 2 \int_{\omega}g_0^+.\overline{b}\, dx_{\alpha},
\end{equation}
where $\overline{b}(x_{\alpha}):=\frac{1}{2}\int_{-1}^{1}b(x_{\alpha},x_3)dx_3$ and $\overline{f}(x_{\alpha}):=\frac{1}{2}\int_{-1}^{1}f(x_{\alpha},x_3)dx_3$. The desired membrane model should thus depend on the average, $\overline{b}$, of $b$ with respect to $x_3$. Once we establish our $\Gamma$-convergence result, we will be in a position to conclude that    $v$ and $\overline{b}$ are truly independent and
that the corresponding model is a Cosserat type membrane model.

To this end, we define, for all $(u,\overline{b};A) \in L^p(\Omega;\mathbb{R}^3) \times  L^p(\omega;\mathbb{R}^3) \times \mathcal{A}(\omega)$,
\begin{equation}\label{funct}
\mathcal{J}_{\varepsilon}(u,\overline{b};A):=\left\{
\begin{array}{ll}
\ds \int_{A \times (-1,1)}W\left(x;D_{\alpha} u\Big|\frac{1}{\varepsilon}D_3 u\right) dx & \text{ if } \left\{
\begin{array}{l}u \in W^{1,p}(A \times (-1,1);\mathbb{R}^3),\\
 \frac{1}{2\varepsilon}\int_{-1}^{1}D_3 u(.,x_3)dx_3=\overline{b},\end{array}\right.\\[5mm]
+\infty & \text{ otherwise},
\end{array}
\right.
\end{equation}
and
\begin{eqnarray}\label{bm:13}
\mathcal{J}_{\{\varepsilon\}}(u,\overline{b};A)  :=  \inf_{\{u_{\varepsilon},\overline{b}_{\varepsilon}\}}\Big\{ \liminf_{\varepsilon \rightarrow 0} \mathcal{J}_{\varepsilon}(u_{\varepsilon},\overline{b}_{\varepsilon};A) : u_{\varepsilon} \rightarrow u \text{ in } L^p(A \times (-1,1);\mathbb{R}^3) \text{ and }\nonumber\\
\overline{b}_{\varepsilon} \rightharpoonup \overline{b} \text{ in } L^p(A;\mathbb{R}^3)\Big\}.
\end{eqnarray}
\begin{rem}\label{w1p2}
Let $(u,\overline{b};A) \in L^p(\Omega;\mathbb{R}^3) \times L^p(\omega;\mathbb{R}^3) \times \mathcal{A}(\omega)$ and suppose that $\mathcal J_{\{\varepsilon\}}(u,\overline{b};A) < +\infty$. Arguing as in Remark \ref{dem2}, we deduce that $u \in W^{1,p}(A;\mathbb{R}^{3})$. Hence, if $u \in L^p(\Omega;\mathbb{R}^3) \setminus W^{1,p}(A;\mathbb{R}^{3})$, then $\mathcal{J}_{\{\varepsilon\}}(u,\overline{b};A)=+\infty$.
\end{rem}
\begin{rem}
Whenever $u \in W^{1,p}(A;\mathbb R^3)$, one has $\mathcal J_{\{\varepsilon\}}(u,\overline{b};A) < +\infty$, which is easily obtained by considering the sequence $\{u_\varepsilon(x_\alpha,x_3):=u(x_\alpha)+\varepsilon x_3 \overline b_\varepsilon(x_\alpha)\}$, where $\overline b_\varepsilon \in \mathcal C^\infty_c(A;\mathbb R^3)$ and $\overline b_\varepsilon \to \overline b$ strongly in $L^p(A;\mathbb R^3)$.
\end{rem}

Theorem 1.2 in  \cite{Bo&Fo&Ma} shows  that, if $W$ is a homogeneous elastic energy density, then  $\mathcal{J}_{\{\varepsilon\}}$ is the $\Gamma(L^p)$-limit of $\mathcal{J}_{\varepsilon}$,  by which we mean, from now onward, the $\Gamma$-limit with respect to, respectively, the strong topology of $L^p(\Omega;\mathbb{R}^3)$, and  the weak topology of $L^p(\omega;\mathbb{R}^3)$. Furthermore, for all $(u,\overline{b};A) \in W^{1,p}(\omega;\mathbb{R}^3) \times L^p(\omega;\mathbb{R}^3) \times \mathcal{A}(\omega)$,
$$\mathcal{J}_{\{\varepsilon\}}(u,\overline{b};A)=2\int_{A}\mathcal{Q}^*W(D_{\alpha}u|\overline{b})dx_{\alpha},$$
where, for all $\overline{F} \in \mathbb{M}^{3 \times 2}$ and $z \in \mathbb{R}^3$,
\begin{eqnarray*}
\mathcal{Q}^*W(\overline{F}|z)  :=  \inf_{L>0,\varphi} \left\{ \frac{1}{2}\int_{Q' \times (-1,1)} W(\overline{F}+D_{\alpha} \varphi|L D_3 \varphi)dx_{\alpha}dx_3 :  \varphi \in W^{1,p}(Q' \times (-1,1);\mathbb{R}^3),\right.\\
\left.\varphi(.,x_3) \; Q'-\text{periodic for a.e. } x_3 \in (-1,1), \frac{L}{2}\int_{Q' \times (-1,1)}D_3\varphi \, dx = z \right\}.
\end{eqnarray*}
We propose to extend this result to the heterogeneous case. We set, for all $\overline{F} \in \mathbb{M}^{3 \times 2}$, $z \in \mathbb{R}^3$ and a.e. $x_0 \in \omega$,
\begin{eqnarray}\label{bm:1}
\mathcal{Q}^*W(x_0;\overline{F}|z) := \inf_{L>0,\varphi} \left\{ \frac{1}{2}\int_{Q' \times (-1,1)} \mathcal{Q}W(x_0,x_3;\overline{F}+D_{\alpha} \varphi|L D_3 \varphi)dx_{\alpha}dx_3 :\right.\nonumber\\ \varphi \in W^{1,p}(Q' \times (-1,1);\mathbb{R}^3), \varphi(.,x_3) \; Q'\text{-periodic for a.e. } x_3 \in (-1,1),\\
\left. \text{and }\frac{L}{2}\int_{Q' \times (-1,1)}D_3\varphi \, dx = z\right\}\nonumber
\end{eqnarray}
where, for a.e. $x \in \Omega$ and all $F \in \mathbb{M}^{3 \times 3}$,
$\mathcal{Q}W(x;.)$, the 3D-quasiconvexification of $W(x;.)$ is defined as
$$\mathcal{Q}W(x;F)=\inf_{\phi \in W^{1,p}_0(Q;\mathbb{R}^3)}\int_Q W(x;F+D\phi(y))dy$$
with $Q:=(0,1)^3$. Since $\mathcal{Q}W(x;.)$ is quasiconvex and satisfies a $p$-growth condition,    for all $F_1,F_2 \in \mathbb{M}^{3 \times 3}$ and for a.e. $x \in \Omega$,
\begin{equation}\label{bm:3}
|\mathcal{Q}W(x;F_1)-\mathcal{Q}W(x;F_2)|\leq \beta (1 +|F_1|^{p-1}+|F_2|^{p-1})|F_1-F_2|
\end{equation}
 (see \cite{Dac1}, Lemma 2.2).  Elementary properties of $\mathcal{Q}^*W$  are summarized in the following proposition:
\begin{prop}
\label{propsup}
i) For all $(\overline{F},z) \in \mathbb{M}^{3 \times 2} \times \mathbb{R}^3$ and a.e. $x_0 \in \omega$,
\begin{equation}\label{bm:2}
0 \leq \mathcal{Q}^*W(x_0;\overline{F}|z) \leq \beta(|\overline{F}|^p+|z|^p+1).
\end{equation}
ii) $\mathcal{Q}^*W$ is a Carath\'eodory function.
\end{prop}
\textit{Proof. } Item i). We take $\varphi(x):=z x_3/L$ as test function in (\ref{bm:1}) an use the $p$-growth condition (\ref{A}).\\
Item ii). It suffices to show that $\mathcal Q^*W(x_0;.)$ is continuous for a.e. $x_0 \in \omega$. Let $\overline{F}_n \rightarrow \overline{F}$ and $z_n \rightarrow z$. We first prove that $\mathcal Q^*W(x_0;.)$ is upper semicontinuous. For any $\delta >0$, set $L>0$ and $\varphi \in W^{1,p}(Q' \times (-1,1);\mathbb{R}^3)$ $Q'$-periodic satisfying $\frac{L}{2}\int_{Q' \times (-1,1)}D_3\varphi dx = z$ such that
$$\mathcal Q^*W(x_0;\overline{F}|z) \leq \frac{1}{2}\int_{Q' \times (-1,1)} \mathcal{Q}W(x_0,x_3;\overline{F}+D_{\alpha} \varphi|L D_3 \varphi)dx_{\alpha}dx_3 \leq \mathcal Q^*W(x_0;\overline{F}|z) + \delta.$$
The sequence $\{\varphi_n(x):=\varphi(x)+x_3(z_n-z)/L\}$ is in $W^{1,p}(Q' \times (-1,1);\mathbb{R}^3)$
and it is $Q'$-periodic. Furthermore, $D_{\alpha}\varphi=D_{\alpha}\varphi_n$ and $\frac{L}{2}\int_{Q' \times (-1,1)}D_3\varphi_n dx = z_n$.
Since
$$\|(\overline{F}+D_{\alpha}\varphi|L D_3\varphi)-(\overline{F}_n+D_{\alpha}\varphi_n|L D_3\varphi_n)\|_{L^{\infty}(Q' \times (-1,1);\mathbb{M}^{3 \times 3})} \leq |(\overline{F}|z)-(\overline{F}_n|z_n)| \rightarrow 0$$
while $\{D \varphi_n\}$ is bounded in $L^p(Q' \times (-1,1);\mathbb{M}^{3 \times 3})$, (\ref{bm:3}) together with H\"older's inequality, yields
\begin{eqnarray*}
& & \limsup_{n \rightarrow +\infty}\mathcal Q^*W(x_0;\overline{F}_n|z_n)-\mathcal Q^*W(x_0;\overline{F}|z)-\delta\\
& \leq &\limsup_{n \rightarrow +\infty} \frac{1}{2}\int_{Q' \times (-1,1)}\left|\mathcal{Q}W(x_0,x_3;\overline{F}_n+D_{\alpha}\varphi_n|L D_3\varphi_n)-\mathcal{Q}W(x_0,x_3;\overline{F}+D_{\alpha}\varphi|L D_3\varphi)\right|dx\\
& \leq & \limsup_{n \rightarrow +\infty} \frac{\beta}{2}\int_{Q' \times (-1,1)}\left(1+|(\overline{F}_n+D_{\alpha}\varphi_n|L D_3\varphi_n)|^{p-1}+|(\overline{F}+D_{\alpha}\varphi|L D_3\varphi)|^{p-1}\right) \\
& & \qquad\qquad \times |(\overline{F}+D_{\alpha}\varphi|L D_3\varphi)-(\overline{F}_n+D_{\alpha}\varphi_n|L D_3\varphi_n)| dx\\
& \leq & \limsup_{n \rightarrow +\infty}C\left(1+\|D\varphi\|^{p-1}_{L^p(Q' \times (-1,1);\mathbb{M}^{3 \times 3})}+\|D\varphi_n\|^{p-1}_{L^p(Q' \times (-1,1);\mathbb{M}^{3 \times 3})}\right) \\
& &\qquad\qquad \times \|(\overline{F}+D_{\alpha}\varphi|L D_3\varphi)-(\overline{F}_n+D_{\alpha}\varphi_n|L D_3\varphi_n)\|_{L^{p}(Q' \times (-1,1);\mathbb{M}^{3 \times 3})}=0
\end{eqnarray*}
Passing to the limit when $\delta \! \searrow \! 0^+$ yields the desired upper semicontinuity. Let us prove now that $\mathcal Q^*W(x_0;.)$ is lower semicontinuous. For every $n \geq 1$, choose $L_n>0$ and $\varphi_n \in W^{1,p}(Q' \times (-1,1);\mathbb{R}^3)$ $Q'$-periodic satisfying $\frac{L_n}{2}\int_{Q' \times (-1,1)}D_3\varphi_n dx = z_n$ such that
$$\frac{1}{2}\int_{Q' \times (-1,1)} \mathcal{Q}W(x_0,x_3;\overline{F}_n+D_{\alpha} \varphi_n|L_n D_3 \varphi_n)dx_{\alpha}dx_3 \leq \mathcal Q^*W(x_0;\overline{F}_n|z_n) + \frac{1}{n}.$$
Set $\tilde \varphi_n(x):=\varphi_n(x)+x_3 (z-z_n)/L_n$, then $\tilde \varphi_n \in W^{1,p}(Q' \times (-1,1);\mathbb{R}^3)$ is $Q'$-periodic and satisfies $\frac{L_n}{2}\int_{Q' \times (-1,1)}D_3\tilde \varphi_n dx = z$. Since
$$\|(\overline{F}+D_{\alpha}\tilde \varphi_n|L_n D_3\tilde \varphi_n)-(\overline{F}_n+D_{\alpha}\varphi_n|L_n D_3\varphi_n)\|_{L^{\infty}(Q' \times (-1,1);\mathbb{M}^{3 \times 3})} \leq |(\overline{F}|z)-(\overline{F}_n|z_n)| \rightarrow 0$$
while, in view of the coercivity condition (\ref{A}), the sequences $\{(D_{\alpha} \varphi_n|L_n D_3 \varphi_n)\}$ and $\{(D_{\alpha} \tilde \varphi_n|L_n D_3 \tilde \varphi_n)\}$ are bounded in $L^p(Q'\times (-1,1);\mathbb M^{3 \times 3})$ uniformly in $n$, (\ref{bm:3}) implies that
\begin{eqnarray*}
\mathcal Q^*W(x_0;\overline F|z) & \leq & \liminf_{n \to +\infty}\frac{1}{2}\int_{Q' \times (-1,1)} \mathcal{Q}W(x_0,x_3;\overline{F}+D_{\alpha} \tilde \varphi_n|L_n D_3 \tilde \varphi_n)dx_{\alpha}dx_3\\
 & \leq & \liminf_{n \to +\infty}\frac{1}{2}\int_{Q' \times (-1,1)} \mathcal{Q}W(x_0,x_3;\overline{F}_n+D_{\alpha} \varphi_n|L_n D_3 \varphi_n)dx_{\alpha}dx_3\\
 & \leq & \liminf_{n \to +\infty}\mathcal Q^*W(x_0;\overline{F}_n|z_n).
\end{eqnarray*}
Thus $ Q^*W(x_0;.)$ is lower semicontinuous and the continuity follows.
\hfill$\Box$\\

We propose to establish the following $\Gamma$-convergence result.
\begin{theo}
\label{t3.1}
For all $A \in \mathcal{A}(\omega)$, $\mathcal{J}_{\varepsilon}(.,.;A)$ $\Gamma(L^p)$-converges to $\mathcal{J}_{\{\varepsilon\}}(.,.;A)$. Further, for all $(u,\overline{b}) \in L^p(\Omega;\mathbb{R}^3) \times L^p(\omega;\mathbb{R}^3)$,
$$\mathcal{J}_{\{\varepsilon\}}(u,\overline{b};A)=\left\{
\begin{array}{ll}
\ds 2\int_{A}\mathcal{Q}^*W(x_{\alpha};D_{\alpha}u|\overline{b})dx_{\alpha} & \text{if } u \in W^{1,p}(A;\mathbb{R}^3),\\
 & \\
+\infty & \text{otherwise},
\end{array}
\right.$$
where $\mathcal{Q}^*W$ is given by (\ref{bm:1}).
\end{theo}

We first note, as in \cite{Bra&Fo&Fr} p.1374, that, if $\mathcal R(\omega)$ is the countable family of all finite unions of open squares in $\omega$ with faces parallel to the axes, centered at rational points and with rational edge lengths, then there exists a subsequence $\{\varepsilon_n\} \subset \{\varepsilon\}$ such that $\mathcal J_{\{\varepsilon_n\}}(.,.;C)$ is, for all $C \in \mathcal R(\omega)$, the $\Gamma(L^p)$-limit of $\mathcal J_{\varepsilon_n}(.,.;C)$.

Then, the analogue of Step 2 in the proof of Theorem 2.5 of \cite{Bra&Fo&Fr} holds, namely
\begin{lem}
\label{l3.0}
For any $A \in \mathcal{A}(\omega)$ and $(u,\overline{b}) \in W^{1,p}(A;\mathbb{R}^3) \times L^p(A;\mathbb{R}^3)$, there exists a sequence $\{u_n\} \subset W^{1,p}(A \times (-1,1);\mathbb{R}^3)$ satisfying
\begin{equation}\label{gl}
\left\{
\begin{array}{l}
u_n \to u \text{ in } L^p(A \times (-1,1);\mathbb{R}^3),\\
\\
\ds \overline{b}_n:= \frac{1}{2\varepsilon_n}\int_{-1}^1 D_3 u_n(.;x_3)dx_3 \rightharpoonup \overline{b} \text{ in } L^p(A ;\mathbb{R}^3),\\
\\
\lim\limits_{n \to +\infty} \mathcal{J}_{\varepsilon_n}\left(u_n,\overline{b}_n;A\right)=\mathcal{J}_{\{\varepsilon_n\}}(u,\overline{b};A).
\end{array}
\right.
\end{equation}
\end{lem}
\noindent\textit{Proof. }The coercivity condition (\ref{A}) implies that whenever $u \in W^{1,p}(C;\mathbb{R}^3)$, we can choose the attainment sequence $\{u_n,\overline{b}_n\}$, so that (\ref{gl}) holds true. Now let us fix $\delta >0$ and choose a subset $C^{\delta}$ of $A$ in $\mathcal{R}(\omega)$ such that $\overline{C^{\delta}} \subset A$ and
$$\int_{A \setminus C^{\delta}}(1+|D_{\alpha}u|^p)dx_{\alpha} < \frac{\delta}{2\beta}.$$
Consider a sequence $\{v^{\delta}_n\} \subset W^{1,p}(C^\delta \times (-1,1);\mathbb{R}^3)$ satisfying
$$\left\{
\begin{array}{l}
v^{\delta}_n \xrightarrow[n \to +\infty]{} u \text{ in } L^p(C^{\delta} \times (-1,1);\mathbb{R}^3),\\
\\
\ds \overline{b}_n^\delta := \frac{1}{2\varepsilon_n}\int_{-1}^1 D_3 v^{\delta}_n(.,x_3)dx_3 \xrightharpoonup[n \to +\infty]{} \overline{b} \text{ in } L^p(C^{\delta} ;\mathbb{R}^3),\\
\\
\lim\limits_{n \to +\infty} \mathcal{J}_{\varepsilon_n}\left(v^{\delta}_n,\overline{b}_n^\delta
;C^{\delta}\right)=\mathcal{J}_{\{\varepsilon_n\}}(u,\overline{b};C^{\delta}).
\end{array}
\right.$$
In view of Lemma 2.2 in \cite{Bo&Fo&Ma} (the proof in our context is identical to that of the homogeneous case), there exists a subsequence of $\{\varepsilon_n\}$ (not relabelled) and a sequence $\{\hat{v}^{\delta}_n\}$ in $W^{1,p}(C^{\delta} \times (-1,1);\mathbb{R}^3)$ satisfying $\hat{v}^{\delta}_n=u$ on a neighborhood of $\partial C^{\delta} \times (-1,1)$ such that
\begin{equation}\label{xxl}
\left\{
\begin{array}{l}
\hat{v}^{\delta}_n \xrightarrow[n \to +\infty]{} u \text{ in } L^p(C^{\delta} \times (-1,1);\mathbb{R}^3),\\
\\
\ds\overline{\hat b}_n^\delta:= \frac{1}{2\varepsilon_n}\int_{-1}^{1}D_3 \hat{v}^{\delta}_n(.;x_3)dx_3 \xrightharpoonup[n \to +\infty]{} \overline{b} \text{ in } L^p(C^{\delta} ;\mathbb{R}^3),\\
\\
\lim\limits_{n \to +\infty} \mathcal{J}_{\varepsilon_n}\left(\hat{v}^{\delta}_n,\overline{\hat b}_n^\delta;C^{\delta}\right)=\mathcal{J}_{\{\varepsilon_n\}}(u,\overline{b};C^{\delta}).
\end{array}
\right.
\end{equation}
We  extend $\hat{v}^{\delta}_n$ as $u$ outside $C^{\delta}$ (and correspondingly extend $\overline{\hat b}_n^\delta$ as $0$). Since $\mathcal{J}_{\{\varepsilon_n\}}(u,\overline{b};.)$ is an increasing set function, we have $\mathcal{J}_{\{\varepsilon_n\}}(u,\overline{b};C^{\delta}) \leq \mathcal{J}_{\{\varepsilon_n\}}(u,\overline{b};A)$ and thus,
\begin{eqnarray*}
& & \limsup_{\delta \to 0^+} \limsup_{n \to +\infty} \; \mathcal{J}_{\varepsilon_n}\left(\hat{v}^{\delta}_n,\overline{\hat b}_n^\delta;A\right)\\
 & \leq & \limsup_{\delta \to 0^+} \left\{ \lim_{n \to +\infty} \mathcal{J}_{\varepsilon_n}\left(\hat{v}^{\delta}_n,\overline{\hat b}_n^\delta;C^{\delta}\right) + 2\beta \int_{A \setminus  C^{\delta}}(1+|D_{\alpha}u|^p)dx_{\alpha} \right\}\\
 & = & \limsup_{\delta \to 0^+}  \mathcal{J}_{\{\varepsilon_n\}}(u,\overline{b};C^{\delta})\\
 & \leq & \mathcal{J}_{\{\varepsilon_n\}}(u,\overline{b};A)\\
 & \leq & \liminf_{\delta \to 0^+} \liminf_{n \to +\infty} \; \mathcal{J}_{\varepsilon_n}\left(\hat{v}^{\delta}_n,\overline{\hat b}_n^\delta;A\right).
\end{eqnarray*}
Remark that (\ref{xxl}), together with coercivity, implies that
$$\|D_\alpha \hat v^{\delta}_n\|_{L^p(A\times(-1,1);\mathbb{M}^{3\times 2})}+\| \overline{\hat b}_n^\delta\|_{L^p(A;\mathbb{R}^{3})} \leq C,$$
independently of $\delta,n$; in particular, $\{\overline{\hat b}_n^\delta\}$ lies in a subset of $L^p(A;\mathbb{R}^3)$,
 which is metrizable for the weak $L^p$-topology. A simple diagonalization lemma (Lemma 7.1 in \cite{Bra&Fo&Fr}) permits to conclude the existence of a decreasing sequence $\{\delta(n)\} \searrow 0^+$ such that the sequence $\{u_n:=\hat{v}^{\delta(n)}_n\}$ satisfies (\ref{gl}). \hfill $\Box$\\

We now recall two results that will be of use in the proof of Lemma \ref{l3.2} below. Their  proof  can be found in \cite{Bo&Fo&Ma} in the homogeneous case and the heterogeneity does not create any additional difficulty.
\begin{prop}\label{propmes}
For any sequence $\{\varepsilon\} \searrow 0^+$, there exists a subsequence $\{\varepsilon_n\}$ such that, for any $(u,\overline{b}) \in W^{1,p}(\omega;\mathbb{R}^3) \times L^p(\omega;\mathbb{R}^3)$, the set function $\mathcal{J}_{\{\varepsilon_n\}}(u,\overline{b};.)$ defined in (\ref{bm:13}) is the trace on $\mathcal{A}(\omega)$ of a Radon measure, which is absolutely continuous with respect to the 2-dimensional Lebesgue measure.
\end{prop}
By virtue of Lemma \ref{l3.0} and Proposition \ref{propmes}, we will assume henceforth that $\{\varepsilon_n\}$ denotes a subsequence of $\{\varepsilon\}$ such that the $\Gamma(L^p)$-limit of $\mathcal{J}_{\varepsilon_n}$ exists, in which case it coincides with $\mathcal{J}_{\{\varepsilon_n\}}$, and such that, for every $(u,\overline{b}) \in W^{1,p}(\omega;\mathbb{R}^3) \times L^p(\omega;\mathbb{R}^3)$, the set function $\mathcal{J}_{\{\varepsilon_n\}}(u,\overline{b};.)$ is the trace on $\mathcal{A}(\omega)$ of a Radon measure, which is absolutely continuous with respect to the 2-dimensional Lebesgue measure..
\begin{prop}\label{propQ}
For all $(u,\overline{b};A) \in W^{1,p}(\omega;\mathbb{R}^3) \times L^p(\omega;\mathbb{R}^3) \times \mathcal{A}(\omega)$, the value  of $\mathcal{J}_{\{\varepsilon\}}(u,\overline{b};A)$ is unchanged if $W$ is replaced by $\mathcal Q W$
in  (\ref{funct}).
\end{prop}
\begin{rem}\label{x3}
If $W$ does not depend on $x_\alpha$, we can show as in \cite{Bo&Fo&Ma} that for all $A \in \mathcal A(\omega)$, $\mathcal{J}_{\varepsilon}(.,.;A)$ $\Gamma(L^p)$-converges to $\mathcal{J}_{\{\varepsilon\}}(.,.;A)$ and
$$\mathcal{J}_{\{\varepsilon\}}(u,\overline{b};A)= 2\int_{A}\mathcal{Q}^*W(D_{\alpha}u|\overline{b})dx_{\alpha},$$
for every $(u,\overline{b}) \in W^{1,p}(A;\mathbb{R}^3) \times L^p(\omega;\mathbb{R}^3)$, where $\mathcal{Q}^*W$ is given by (\ref{bm:1}).
\end{rem}
Theorem \ref{t3.1} is a direct consequence of the following two lemmata.
\begin{lem}
\label{l3.1}
For all $A \in \mathcal{A}(\omega)$ with $A$ Lipschitz and for all $(u,\overline{b}) \in W^{1,p}(A;\mathbb{R}^3) \times L^p(A;\mathbb{R}^3)$,
$$\mathcal{J}_{\{\varepsilon_n\}}(u,\overline{b};A) \geq 2\int_{A}\mathcal{Q}^*W(x_{\alpha};D_{\alpha}u|\overline{b})dx_{\alpha}.$$
\end{lem}
\begin{lem}
\label{l3.2}
For all $A \in \mathcal{A}(\omega)$ with $A$ Lipschitz and for all $(u,\overline{b}) \in W^{1,p}(A;\mathbb{R}^3) \times L^p(A;\mathbb{R}^3)$,
$$\mathcal{J}_{\{\varepsilon_n\}}(u,\overline{b};A) \leq  2\int_{A}\mathcal{Q}^*W(x_{\alpha};D_{\alpha}u|\overline{b})dx_{\alpha}.$$
\end{lem}
\noindent\textit{Proof of Lemma \ref{l3.1}. }Let $(u,\overline{b}) \in W^{1,p}(A;\mathbb{R}^{3}) \times L^p(A;\mathbb{R}^{3})$. According to the equi-integrability Theorem (Theorem 1.1 in \cite{Bo&Fo}) together with Lemma \ref{l3.0} , there exists a subsequence of $\{\varepsilon_n\}$ (not relabelled) and a sequence $\{u_n\} \subset W^{1,p}(A \times (-1,1);\mathbb R^3)$ such that
$$\left\{ \begin{array}{l}
u_n \to u \text{ in }L^p(A \times (-1,1);\mathbb R^3),\\\\
\ds \frac{1}{2\varepsilon_n}\int_{-1}^1 D_3 u_n(.,x_3)dx_3 \rightharpoonup \overline b \text{ in }L^p(A ;\mathbb R^3),\\\\
\left\{ \left| \left( D_\alpha u_n \Big| \frac{1}{\varepsilon_n} D_3 u_n\right) \right|^p \right\} \text{ is equi-integrable},\\\\
\ds \mathcal J_{\{\varepsilon_n\}}(u,\overline b;A) = \lim_{n \to +\infty} \int_{A \times (-1,1)}W\left(x_\alpha,x_3;D_\alpha u_n \Big| \frac{1}{\varepsilon_n} D_3 u_n\right)dx_\alpha dx_3.
\end{array} \right.$$
We argue as in the proof of Lemma \ref{lem>} with $F_n(x):= \left(D_\alpha u_n(x) \big| \frac{1}{\varepsilon_n} D_3 u_n(x)\right)$ and we obtain, since $W \geq \mathcal Q W$
\begin{eqnarray*}
\mathcal J_{\{\varepsilon_n\}}(u,\overline b;A) & \geq & \limsup_{h \to +\infty} \sum_{i \in I_h} \frac{1}{\mathcal L^2(A_{i,h})} \int_{A_{i,h}} \left\{ \liminf_{n \to +\infty} \int_{A_{i,h} \times (-1,1)}  W(y_\alpha,x_3;F_n(x_\alpha,x_3))dx_\alpha dx_3 \right\} dy_\alpha\\
 & \geq & \limsup_{h \to +\infty} \sum_{i \in I_h} \frac{1}{\mathcal L^2(A_{i,h})} \int_{A_{i,h}} \left\{ \liminf_{n \to +\infty} \int_{A_{i,h} \times (-1,1)} \mathcal Q W(y_\alpha,x_3;F_n(x_\alpha,x_3))dx_\alpha dx_3 \right\} dy_\alpha,
\end{eqnarray*}
where $\{A_{i,h}\}_{i \in I_h}$ denotes a finite family of disjoint open sets with diameter less than $1/h$, such that $\mathcal L^2(A \setminus \cup_{i \in I_h}A_{i,h})=0$ and ${\rm Card }I_h =O(h^2)$. Applying, for a.e. $y_\alpha \in A_{i,h}$, Remark \ref{x3} to the Carath\'eodory function $(x_3;F) \mapsto W(y_\alpha,x_3;F)$, we get
$$\mathcal J_{\{\varepsilon_n\}}(u,\overline b;A) \geq \limsup_{h \to +\infty} \sum_{i \in I_h} \frac{1}{\mathcal L^2(A_{i,h})} \int_{A_{i,h}} \left\{ 2 \int_{A_{i,h}}\mathcal Q^*W(y_\alpha;D_\alpha u(x_\alpha)|\overline b(x_\alpha))dx_\alpha \right\} dy_\alpha.$$
By Proposition \ref{propsup} (ii), $\mathcal Q^*W$ is a Carath\'eodory integrand, thus Scorza-Dragoni's Theorem implies the existence, for any $\eta>0$, of a compact set $C_\eta \subset A$, such that
\begin{equation}\label{ceta}
\mathcal L^2(A \setminus C_\eta)<\eta,
\end{equation}
and $\mathcal Q^*W$ is continuous on $C_\eta \times \mathbb M^{3 \times 2}$. Let $S_\lambda:=\{ x_\alpha \in A : |(D_\alpha u(x_\alpha)|\overline b(x_\alpha))| \leq \lambda\}$, thanks to Chebyshev's inequality
\begin{equation}\label{slambda}
\mathcal L^2(A \setminus S_\lambda) <  \frac{C}{\lambda^p}.
\end{equation}
Consequently
\begin{eqnarray*}
&\ds \mathcal J_{\{\varepsilon_n\}}(u,\overline b;A) \geq\\
&\ds \limsup_{\lambda \to +\infty} \limsup_{\eta \to 0} \limsup_{h \to +\infty} \sum_{i \in I_h} \frac{1}{\mathcal L^2(A_{i,h})} \int_{A_{i,h} \cap C_\eta} \left\{ 2 \int_{A_{i,h} \cap C_\eta \cap S_\lambda}\mathcal Q^*W(y_\alpha;D_\alpha u(x_\alpha)|\overline b(x_\alpha))dx_\alpha \right\} dy_\alpha.
\end{eqnarray*}
Since $\mathcal Q^*W$ is continuous on $C_\eta \times \mathbb M^{3 \times 2}$, it is uniformly continuous on $C_\eta \times \overline B(0,\lambda)$ thus there exists a increasing and continuous function $\omega_{\eta,\lambda} :\mathbb R^+ \longrightarrow \mathbb R^+$ satisfying $\omega_{\eta,\lambda}(0)=0$ and for every $y_\alpha \in A_{i,h} \cap C_\eta$ and every $x_\alpha \in A_{i,h} \cap C_\eta \cap S_\lambda$,
$$\big|\mathcal Q^*W(y_\alpha;D_\alpha u(x_\alpha)|\overline b(x_\alpha))-\mathcal Q^*W(x_\alpha;D_\alpha u(x_\alpha)|\overline b(x_\alpha))\big| \leq \omega_{\eta,\lambda}(|x_\alpha-y_\alpha|) \leq \omega_{\eta,\lambda}(1/h).$$
Using the fact that $ \mathcal L^2(A_{i,h}){\rm Card }I_h \leq C$, we get
\begin{eqnarray*}
&\ds \sum_{i \in I_h}  \frac{1}{\mathcal L^2(A_{i,h})} \int_{A_{i,h} \cap C_\eta} \int_{A_{i,h} \cap C_\eta \cap S_\lambda} \big|\mathcal Q^*W(y_\alpha;D_\alpha u(x_\alpha)|\overline b(x_\alpha))-\mathcal Q^*W(x_\alpha;D_\alpha u(x_\alpha)|\overline b(x_\alpha))\big| dx_\alpha dy_\alpha\\
&\ds \leq C\omega_{\eta,\lambda}(1/h) \xrightarrow[h \to +\infty]{}0.
\end{eqnarray*}
Therefore,
$$ \mathcal J_{\{\varepsilon_n\}}(u,\overline b;A) \geq  \limsup_{\lambda \to +\infty} \limsup_{\eta \to 0} \limsup_{h \to +\infty} \sum_{i \in I_h}   \frac{2 \mathcal L^2(A_{i,h} \cap C_\eta)}{\mathcal L^2(A_{i,h})} \int_{A_{i,h} \cap C_\eta \cap S_\lambda}\mathcal Q^*W(x_\alpha;D_\alpha u(x_\alpha)|\overline b(x_\alpha))dx_\alpha.$$
By virtue of the $p$-growth condition (\ref{bm:2}) together with (\ref{ceta}), we get
\begin{eqnarray*}
\sum_{i \in I_h} \frac{\mathcal L^2(A_{i,h} \setminus C_\eta)}{\mathcal L^2(A_{i,h})}  \int_{A_{i,h} \cap C_\eta \cap S_\lambda}\mathcal Q^*W(x_\alpha;D_\alpha u(x_\alpha)|\overline b(x_\alpha))dx_\alpha & \leq & \beta(1+\lambda^p)\sum_{i \in I_h} \mathcal L^2(A_{i,h} \setminus C_\eta)\\
 & = & \beta(1+\lambda^p)\mathcal L^2(A \setminus C_\eta)\\
 & < & \beta(1+\lambda^p)\eta \xrightarrow[\eta \to 0]{} 0.
\end{eqnarray*}
Thus, (\ref{ceta}) and (\ref{slambda}) yield
\begin{eqnarray*}
\mathcal J_{\{\varepsilon_n\}}(u,\overline b;A) & \geq & 2 \limsup_{\lambda \to +\infty} \limsup_{\eta \to 0} \int_{A \cap C_\eta \cap S_\lambda}\mathcal Q^*W(x_\alpha;D_\alpha u(x_\alpha)|\overline b(x_\alpha))dx_\alpha\\
&=&  2 \int_{A}\mathcal Q^*W(x_\alpha;D_\alpha u(x_\alpha)|\overline b(x_\alpha))dx_\alpha.
\end{eqnarray*}
\hfill$\Box$
\\
\\
\noindent\textit{Proof of Lemma \ref{l3.2}. }The proof is divided into three steps. First, we address the case where $u$ is affine and $\overline{b}$ is constant; then, that where $u$ is piecewise affine and continuous, and $\overline{b}$ piecewise constant. Finally, we address  the general case.\\
\textit{Step 1. }Let $A  \in \mathcal A(\omega)$, we assume that
$$\left\{
\begin{array}{ll}
u(x_{\alpha})=\overline{F}x_{\alpha}+c, &  (\overline{F},c) \in \mathbb{M}^{3 \times 2} \times \mathbb R^3,\\
\overline{b}(x_{\alpha})=z, & z \in  \mathbb{R}^3.
\end{array}
\right.$$
Thanks to the density of $W^{1,\infty}(Q'\times (-1,1);\mathbb{R}^3)$ into $W^{1,p}(Q'\times (-1,1);\mathbb{R}^3)$ and to the $p$-growth hypothesis (\ref{A}), for any $k \geq 1$, there exists $L_k>0$ and $\varphi_k \in W^{1,\infty}(Q'\times (-1,1);\mathbb{R}^3)$ $Q'$-periodic satisfying $\frac{L_k}{2}\int_{Q'\times(-1,1)}D_3\varphi_k dx=z$ and such that
$$Z_k(x_0;\overline F|z):=\frac{1}{2}\int_{Q'\times(-1,1)}\mathcal{Q} W\left( x_0,x_3;\overline{F}+D_{\alpha}\varphi_k |L_k D_3\varphi_k \right)dx_{\alpha}dx_3 \leq \mathcal{Q}^*W(x_0;\overline{F}|z)+\frac{1}{k}.$$
We extend $\varphi_k$ to $\mathbb R^2 \times (-1,1)$ by $Q'$-periodicity. Choose $r>0$ small enough so that $Q'(x_0,r) \subset A$ where $x_0$ is a Lebesgue point of the Radon-Nikodym derivative of $\mathcal{J}_{\{\varepsilon_n\}}(u,\overline{b};.)$ with respect to the 2-dimensional Lebesgue measure and of $Z_k(.;\overline F|z)$ for all $k \geq 1$. Fix $k$ and set
$$u_n^k(x):=\overline{F}x_{\alpha}+c+L_k\varepsilon_n\varphi_k\left(\frac{x_{\alpha}}{L_k\varepsilon_n},x_3\right).$$
Then,
$$u_n^k \xrightarrow[n \to +\infty]{} u \text{ in } L^p(Q'(x_0,r) \times (-1,1);\mathbb{R}^3),$$
and by virtue of  Riemann-Lebesgue's Lemma,
$$\frac{1}{2\varepsilon_n}\int_{-1}^{1}D_3u_n^k dx_3=\frac{L_k}{2}\int_{-1}^{1}D_3\varphi_k\left(\!\frac{x_{\alpha}}{L_k\varepsilon_n},x_3\!\right)dx_3 \xrightharpoonup[n \to +\infty]{ L^p(Q'(x_0,r);\mathbb{R}^3)} \frac{L_k}{2}\int_{Q'\times(-1,1)}D_3\varphi_k dx\!=\!\overline b.$$
So $\{u_n^k\}$ is admissible for $\mathcal{J}_{\{\varepsilon_n\}}(u,\overline{b};Q'(x_0,r))$ and, thanks to Proposition \ref{propQ},
$$\mathcal{J}_{\{\varepsilon_n\}}(u,\overline{b};Q'(x_0,r)) \leq \liminf_{n \rightarrow +\infty} \int_{Q'(x_0,r) \times (-1,1)}\mathcal{Q}W\left(x_{\alpha},x_3;D_{\alpha}u_n^k \Big|\frac{1}{\varepsilon_n}D_3u_n^k\right)dx_{\alpha}dx_3.$$
Using an argument similar to that in the proof of Lemma \ref{lem<}, with $\mathcal Q W$ instead of $W$, we get
$$\frac{d\mathcal{J}_{\{\varepsilon_n\}}(u,\overline{b};.)}{d\mathcal{L}^2}(x_0) \leq 2\mathcal{Q}^*W(x_0;\overline{F}|z).$$
Thus, integration over $A$ yields
$$\mathcal{J}_{\{\varepsilon_n\}}(u,\overline{b};A) \leq 2\int_{A}\mathcal{Q}^*W(x_{\alpha};\overline{F}|z)dx_{\alpha}.$$
\\
\textit{Step 2. }Assume that $u$ is continuous and piecewise affine and $\overline{b}$ is piecewise constant on $A$. There exists a partition  $A_1,...,A_N$ of $A$ such that $u(x_{\alpha})=\overline{F}_ix_{\alpha}+c_i$ and $\overline{b}(x_{\alpha})=z_i$ on $A_i$, for all $i=1,...,N$. Thanks to step 1, for all $i=1,...,N$, we have
$$\mathcal{J}_{\{\varepsilon_n\}}(\overline{F}_i x_{\alpha}+c_i,z_i;A_i) \leq 2\int_{A_i}\mathcal{Q}^*W(x_{\alpha};\overline{F}_i|z_i)dx_{\alpha}.$$
In view of Proposition \ref{propmes}, $\mathcal{J}_{\{\varepsilon_n\}}(u,\overline{b};.)$ is a measure and we thus get
\begin{eqnarray*}
\mathcal{J}_{\{\varepsilon_n\}}(u,\overline{b};A) & = & \sum_{i=1}^{N}\mathcal{J}_{\{\varepsilon_n\}}(\overline{F}_i x_{\alpha}+c_i,z_i;A_i) \\
 & \leq & 2\sum_{i=1}^{N}\int_{A_i}\mathcal{Q}^*W(x_{\alpha};\overline{F}_i|z_i)dx_{\alpha}\\
 & = &2\int_{A}\mathcal{Q}^*W(x_{\alpha};D_{\alpha}u|\overline{b})dx_{\alpha}.
\end{eqnarray*}
\\
\textit{Step 3. }Consider $A \in \mathcal{A}(\omega)$ with  $A$ Lipschitz and $u \in W^{1,p}(A;\mathbb{R}^3)$, $\overline{b} \in L^p(A;\mathbb{R}^3)$. There exists a sequence $\{u_n\}$ of continuous and piecewise affine functions in $W^{1,p}(A;\mathbb{R}^3)$ and a sequence $\{\overline{b}_n\}$ of piecewise constant functions in $L^p(A;\mathbb{R}^3)$ such that $u_n \rightarrow u$ in $W^{1,p}(A;\mathbb{R}^3)$ and $\overline{b}_n \rightarrow \overline{b}$ in $L^p(A;\mathbb{R}^3)$. Since $\mathcal{J}_{\{\varepsilon_n\}}(.,.;A)$ is lower semicontinuous, we get, in
view of the previous step,
\begin{equation}\label{bm:32}
\mathcal{J}_{\{\varepsilon_n\}}(u,\overline{b};A) \leq \liminf_{n  \rightarrow +\infty}\mathcal{J}_{\{\varepsilon_n\}}(u_n,\overline{b}_n;A) \leq \liminf_{n  \rightarrow +\infty} 2\int_{A}\mathcal{Q}^*W(x_{\alpha};D_{\alpha}u_n|\overline{b}_n)dx_{\alpha}.
\end{equation}
By Proposition \ref{propsup} and Lebesgue's Dominated Convergence Theorem,
\begin{equation}\label{bm:33}
\lim_{n  \rightarrow +\infty} \int_{A}\mathcal{Q}^*W(x_{\alpha};D_{\alpha}u_n|\overline{b}_n)dx_{\alpha} = \int_{A}\mathcal{Q}^*W(x_{\alpha};D_{\alpha}u|\overline{b})dx_{\alpha}.
\end{equation}
Thus (\ref{bm:32}) and (\ref{bm:33}) yield
$$\mathcal{J}_{\{\varepsilon_n\}}(u,\overline{b};A)  \leq  2\int_{A}\mathcal{Q}^*W(x_{\alpha};D_{\alpha}u|\overline{b})dx_{\alpha}.$$
\hfill$\Box$
\\
\\
\textit{Proof of Theorem \ref{t3.1}. }The two previous lemmata demonstrate that, provided $A \in \mathcal{A}(\omega)$ is Lipschitz, then, for all $(u,\overline{b}) \in W^{1,p}(A;\mathbb{R}^3) \times L^p(A;\mathbb{R}^3)$, $\mathcal{J}_{\{\varepsilon_n\}}(u,\overline{b};A)$ does not depend upon the choice of sequence $\{\varepsilon_n\}$. Thus, in light of Proposition 7.11 in \cite{Bra&De}, the whole sequence $\mathcal{J}_{\varepsilon}(u,\overline{b};A)$ $\Gamma(L^p)$-converges to $\mathcal{J}_{\{\varepsilon\}}(u,\overline{b};A)$ and
$$\mathcal{J}_{\{\varepsilon\}}(u,\overline{b};A)=2\int_{A} \mathcal{Q}^*W(x_{\alpha};D_{\alpha}u;\overline{b})dx_{\alpha}.$$
Whenever $A \in \mathcal{A}(\omega)$ is an arbitrary open set, we define the nested sequence of Lipschitz open subsets $A_k:=\{x_\alpha \in A : {\rm dist}(x_\alpha,\partial A) > 1/k\}$ of $A$, so that  $\overline{A_k} \subset A$ and $\cup_{k \geq 1}A_k=A$. Consider the sequence $\{u_k\}$ in $W^{1,p}(A;\mathbb{R}^3)$ such that $u_k=u$ on $A_k$. Since $\mathcal{J}_{\{\varepsilon\}}(.,.;A_k)$ is local and $A_k$ is Lipschitz,
$$\mathcal{J}_{\{\varepsilon\}}(u,\overline{b};A_k)=\mathcal{J}_{\{\varepsilon\}}(u_k,\overline{b};A_k)=2\int_{A_k}\mathcal{Q}^*W(x_{\alpha};D_{\alpha}u_k|\overline{b})dx_{\alpha}=2\int_{A_k}\mathcal{Q}^*W(x_{\alpha};D_{\alpha}u|\overline{b})dx_{\alpha}.$$
But $\mathcal{J}_{\{\varepsilon\}}(u,\overline{b};.)$ is a measure, thus, letting $k \nearrow +\infty$,
$$\mathcal{J}_{\{\varepsilon\}}(u,\overline{b};A)=2\int_{A}\mathcal{Q}^*W (x_{\alpha};D_{\alpha}u|\overline{b})dx_{\alpha}.$$
Then Remark \ref{w1p2} completes the proof of Theorem \ref{t3.1}.\hfill $\Box$
\begin{rem}
If $W$ does not depend upon $x$,  Proposition 1.1 (iii) of \cite{Bo&Fo&Ma} states that
\begin{eqnarray*}
\mathcal{Q}^*W(\overline{F}|z)  =  \inf_{L>0,\varphi} \Big\{ \frac{1}{2}\int_{Q' \times (-1,1)} W(\overline{F}+D_{\alpha} \varphi|L D_3 \varphi)dx_{\alpha}dx_3 :  \varphi \in W^{1,p}(Q' \times (-1,1);\mathbb{R}^3),\\
  \varphi(.,x_3) \; Q'\text{-periodic for a.e. } \; x_3 \in (-1,1), \frac{L}{2}\int_{Q' \times (-1,1)}D_3\varphi dx = z \Big\}\\
 =  \inf_{L>0,\varphi} \Big\{ \frac{1}{2}\int_{Q' \times (-1,1)} \mathcal{Q}W(\overline{F}+D_{\alpha} \varphi|L D_3 \varphi)dx_{\alpha}dx_3 :  \varphi \in W^{1,p}(Q' \times (-1,1);\mathbb{R}^3),\\
  \varphi(.,x_3) \; Q'\text{-periodic for a.e. } \; x_3 \in (-1,1), \frac{L}{2}\int_{Q' \times (-1,1)}D_3\varphi dx = z \Big\}.
\end{eqnarray*}
In other words, the result of \cite{Bo&Fo&Ma} is recovered by Theorem \ref{t3.1}.
\end{rem}
\begin{rem}
\label{r3.2}
Since $\mathcal{Q}^*W$ is the integrand of the $\Gamma(L^p)$-limit of $\mathcal{J}_{\varepsilon}$, which satisfies a $p$-coercivity condition (see (\ref{A})),  for all $\overline{F} \in \mathbb{M}^{3 \times 2}$, for all $z \in \mathbb{R}^3$ and for a.e. $x_0 \in \omega$,
\begin{equation}\label{bm:35}
\beta '(|\overline{F}|^p +|z|^p) \leq \mathcal{Q}^*W(x_0;\overline{F}|z).
\end{equation}
\end{rem}
\begin{rem}
\label{r3.3}
 Theorem \ref{t3.1} implies that the functional
$$(u,\overline{b}) \mapsto \int_{\omega}\mathcal{Q}^*W(x_{\alpha};D_{\alpha}u|\overline{b})dx_{\alpha}$$
is sequentially weakly lower semicontinuous on $W^{1,p}(\omega;\mathbb{R}^3) \times L^p(\omega;\mathbb{R}^3)$. Therefore,  $\mathcal{Q}^*W(x_0;.|z)$ is quasiconvex and $\mathcal{Q}^*W(x_0;\overline{F}|.)$ is convex. Thanks to the $p$-growth condition (\ref{bm:2}), $\mathcal{Q}^*W(x_0;.|.)$ is locally Lipschitz, because it is separately convex (see Theorem 2.3 in \cite{Dac1})
\end{rem}

\section{Classical membrane model obtained as a zero bending moment density}

This section investigates  the coherence of our results. In the absence of a bending moment density ($g_0=0$), we show
below that Theorem \ref{t3.1} boils down to Theorem \ref{t2.1}. We first give another form of the energy density $\underline{W}$ similar to the definition  of $\mathcal{Q}^*W$(see (\ref{bm:1})). Specifically,
\begin{prop}\label{p1}
For all $\overline{F} \in \mathbb{M}^{3 \times 2}$ and for a.e. $x_0 \in \omega$,
\begin{eqnarray*}
\underline{W}(x_0;\overline{F}) & = & \inf_{L,\varphi} \Big\{ \frac{1}{2}\int_{Q' \times (-1,1)} \mathcal{Q}W(x_0,x_3;\overline{F}+D_{\alpha} \varphi|L D_3 \varphi)dx_{\alpha}dx_3 : \\
 & & L>0, \varphi \in W^{1,p}(Q' \times (-1,1);\mathbb{R}^3), \varphi(.,x_3) \; Q'-\text{periodic for a.e. } \; x_3 \in (-1,1) \Big\}.\nonumber
\end{eqnarray*}
\end{prop}
\textit{Proof. } Set
\begin{eqnarray*}
W^*(x_0;\overline{F}) & := & \inf_{L,\varphi} \Big\{ \frac{1}{2}\int_{Q' \times (-1,1)} \mathcal{Q}W(x_0,x_3;\overline{F}+D_{\alpha} \varphi|L D_3 \varphi)dx_{\alpha}dx_3 : \\
 & & L>0, \varphi \in W^{1,p}(Q' \times (-1,1);\mathbb{R}^3), \varphi(.,x_3) \; Q'-\text{periodic for a.e. } \; x_3 \in (-1,1) \Big\}.
\end{eqnarray*}
It is obvious that $\underline{W}(x_0;\overline{F}) \geq W^*(x_0;\overline{F})$. Conversely, for any $\delta>0$, consider $L>0$ and $\varphi \in W^{1,p}(Q' \times (-1,1);\mathbb{R}^3)$ $Q'$-periodic, such that
$$\frac{1}{2}\int_{Q' \times (-1,1)} \mathcal{Q}W(x_0,x_3;\overline{F}+D_{\alpha} \varphi|L D_3 \varphi)dx_{\alpha}dx_3 \leq W^*(x_0;\overline{F}) + \delta.$$
We extend $\varphi$ by $Q'$-periodicity and we set $\varphi_n(x_{\alpha},x_3):=\frac{1}{n}\varphi(nx_{\alpha},x_3)$. Then, Riemann-Lebesgue's Lemma applied to $\int_{-1}^1 W(x_0,x_3;\overline{F}+D_{\alpha} \varphi(.,x_3)|L D_3 \varphi(.,x_3))dx_3$, implies that
\begin{equation}\label{r3}
\lim_{n \rightarrow +\infty}\frac{1}{2}\int_{Q' \times (-1,1)} \mathcal{Q}W(x_0,x_3;\overline{F}+D_{\alpha} \varphi_n|Ln D_3 \varphi_n)dx_{\alpha}dx_3 \leq W^*(x_0;\overline{F}) + \delta.
\end{equation}
For   fixed  $n$,  the relaxation theorem of \cite{Ac&Fu} (see Statement III.7 in \cite{Ac&Fu})-- applied to $Q'\times \left(-\frac{1}{Ln},\frac{1}{Ln}\right)$ and to $\psi_n(x_\alpha,x_3):= \varphi_n(x_\alpha,Lnx_3)$ -- yields a sequence $$\varphi_{n,k} \xrightharpoonup[k \rightarrow +\infty]{} \varphi_n \text{ in } W^{1,p}(Q' \times (-1,1);\mathbb{R}^3)$$ such that,
\begin{eqnarray}\label{r31}
&&\int_{Q' \times (-1,1)} \!\!\mathcal{Q}W(x_0,x_3;\overline{F}+D_{\alpha} \varphi_n|Ln D_3 \varphi_n)dx_{\alpha}dx_3 \nonumber\\
&= &\lim_{k \rightarrow +\infty}\int_{Q' \times (-1,1)} \!\!W(x_0,x_3;\overline{F}+D_{\alpha} \varphi_{n,k}|Ln D_3 \varphi_{n,k})dx_{\alpha}dx_3.
\end{eqnarray}
Thus (\ref{r3}) together with (\ref{r31}) give
$$\lim_{n \rightarrow +\infty}\lim_{k \rightarrow +\infty}\frac{1}{2}\int_{Q' \times (-1,1)} W(x_0,x_3;\overline{F}+D_{\alpha} \varphi_{n,k}|Ln D_3 \varphi_{n,k})dx_{\alpha}dx_3 \leq W^*(x_0;\overline{F}) + \delta.$$
Furthermore, we have,
$$\lim_{n \rightarrow +\infty}\lim_{k \rightarrow +\infty}\|\varphi_{n,k}\|_{L^p(Q' \times (-1,1);\mathbb{R}^3)}=0.$$
By a standard diagonalization process, we can find an increasing sequence $\{k(n)\}$, with $k(n)\stackrel{n}{\nearrow} +\infty$ such that, if we set $\phi_n:=\varphi_{n,k(n)}$,
\begin{equation}\label{r34}
\lim_{n \rightarrow +\infty}\frac{1}{2}\int_{Q' \times (-1,1)} W(x_0,x_3;\overline{F}+D_{\alpha} \phi_n|Ln D_3 \phi_n)dx_{\alpha}dx_3 \leq W^*(x_0;\overline{F}) + \delta,
\end{equation}
and $\phi_n \rightarrow 0$ in $L^p(Q' \times (-1,1);\mathbb{R}^3)$. By virtue of the coercivity hypothesis (\ref{A}),
$$\|( D_{\alpha} \phi_n|Ln D_3 \phi_n)\|_{L^p(Q' \times (-1,1);\mathbb{M}^{3 \times 3})} \leq C.$$
We define the following sequence of non negative bounded Radon measures
$$\lambda_n:=\left(1+|( D_{\alpha} \phi_n|Ln D_3 \phi_n)|^p\right)\chi_{Q' \times (-1,1)} \mathcal{L}^3. $$
The coercive character (\ref{A}) of $W$ permits to assert that, up to a subsequence (not relabelled), there exists a non negative bounded Radon measure $\lambda$ such that
$$\lambda_n \stackrel{*}{\rightharpoonup} \lambda \text{ in } \mathcal{M}_b(\mathbb{R}^3).$$
We  cut $\phi_n$ near the lateral boundary to obtain a sequence which vanishes on $\partial Q' \times (-1,1)$. Let $\theta_k \in \mathcal{C}^{\infty}_c(Q')$ a cut-off function defined by
\begin{equation}\label{r4}
\left\{
\begin{array}{rcl}
\theta_k(x_{\alpha}) & := & \left\{
\begin{array}{lcl}
1 & \text{if} & x_{\alpha} \in Q'(0,1-1/k),\\
0 & \text{if} & x_{\alpha} \notin Q'(0,1-1/(k+1)),
\end{array}\right.\\
\| D_{\alpha}\theta_k\|_{L^{\infty}(Q')} & \leq & Ck^2,
\end{array}
\right.
\end{equation}
We set $\phi_n^k:=\theta_k\phi_n$, since $\phi_n^k=0$ on $\partial Q' \times (-1,1)$, (\ref{r34}) together with (\ref{r4}) yields
\begin{eqnarray}
W^*(x_0;\overline{F}) & \geq & \liminf_{k \rightarrow +\infty} \liminf_{n \rightarrow +\infty} \frac{1}{2}\int_{Q'(0,1-\frac{1}{k}) \times (-1,1)} W(x_0,x_3;\overline{F}+D_{\alpha} \phi_n^k|Ln D_3 \phi_n^k)dx_{\alpha}dx_3 - \delta \nonumber\\
 & \geq & \liminf_{k \rightarrow +\infty} \liminf_{n \rightarrow +\infty} \frac{1}{2}\int_{Q' \times (-1,1)} W(x_0,x_3;\overline{F}+D_{\alpha} \phi_n^k|Ln D_3 \phi_n^k)dx_{\alpha}dx_3\nonumber\\
 & - & \limsup_{k \rightarrow +\infty} \limsup_{n \rightarrow +\infty}\! \frac{1}{2}\!\int_{(Q'(0,1-\frac{1}{k+1})\setminus Q'(0,1-\frac{1}{k})) \times (-1,1)} \!\!\!\!\!\!\!\!W(x_0,x_3;\overline{F}\!\!+\!\!D_{\alpha} \phi_n^k|Ln D_3 \phi_n^k)dx_{\alpha}dx_3\nonumber\\
 & - & \beta(1+|\overline{F}|^p)\limsup_{k \rightarrow +\infty} \mathcal{L}^2\left(Q'\setminus Q'\left(0,1-\frac{1}{k+1}\right)\right) - \delta \nonumber\\
 & \geq & \underline{W}(x_0;\overline{F})-\limsup_{k \rightarrow +\infty} \limsup_{n \rightarrow +\infty} \left\{C \lambda_n\!\!\left( \left(\!Q'\left(0,1\!-\!\frac{1}{k+1}\right)\! \setminus \! Q'\left(0,1\!-\!\frac{1}{k}\right)\right) \! \times \! (-1,1)\right)\right.\nonumber\\
 & + & \left. C' k^{2p} \int_{Q' \times (-1,1)}|\phi_n|^pdx \right\} -\delta. \label{r5}
\end{eqnarray}
Since $Q'(0,1-1/k)$ is an increasing sequence of open sets, the union of which is $Q'$, we get
\begin{eqnarray*}
&& \limsup_{k \rightarrow +\infty}\limsup_{n \rightarrow +\infty} \; \lambda_n\left( \left(Q'\left(0,1-\frac{1}{k+1}\right)\setminus Q'\left(0,1-\frac{1}{k}\right)\right) \times (-1,1)\right)\\
&& \leq  \limsup_{k \rightarrow +\infty} \; \lambda \left( \overline{ \left(Q'\left(0,1-\frac{1}{k+1}\right) \setminus Q'\left(0,1- \frac{1}{k}\right)\right) \times (-1,1) }\right)\\
&& \leq  \limsup_{k \rightarrow +\infty} \; \lambda \left(  \left( Q' \setminus Q'\left(0,1-\frac{1}{k-1}\right)\right) \times [-1,1] \right) =0.
\end{eqnarray*}
Using the fact that $\phi_n \rightarrow 0$ in $L^p(Q' \times (-1,1);\mathbb{R}^3)$ and letting $\delta$ tend to $0$ in (\ref{r5}),
we finally get
$$W^*(x_0;\overline{F}) \geq \underline{W}(x_0;\overline{F}).$$\hfill$\Box$
\\
\\
Now that $\underline{W}$ and $\mathcal{Q}^*W$ are expressed in near identical manner, Remarks \ref{r3.2} and \ref{r3.3} immediately imply that for all $\overline{F} \in \mathbb{M}^{3 \times 2}$ and for a.e. $x_0 \in \omega$, there exists $b_0 \in \mathbb{R}^3$ such that
$$\underline{W}(x_0;\overline{F}) = \min_{z \in \mathbb{R}^3}\mathcal{Q}^*W(x_0;\overline{F}|z)=\mathcal{Q}^*W(x_0;\overline{F}|b_0).$$
In the absence of bending moments, the linear form $L$ given by (\ref{bm:0}) does not depend upon $\overline{b}$ and we may perform explicitly the minimum in $\overline{b}$ in the limit minimization problem. For $u \in W^{1,p}(\omega;\mathbb{R}^3)$, a classical measurability selection criterion (see \cite{Ek&Te}, Chapter VIII, Theorem 1.2),  together with the coercivity condition (\ref{bm:35}), implies the existence of $\overline{b}_0 \in L^p(\omega;\mathbb{R}^3)$ such that for a.e. $x_0 \in \omega$,
$$\underline{W}(x_0;D_{\alpha}u(x_0))=\min_{z \in \mathbb{R}^3}\mathcal{Q}^*W(x_0;D_{\alpha}u(x_0)|z)=\mathcal{Q}^*W(x_0;D_{\alpha}u(x_0)|\overline{b}_0(x_0)).$$
Thus,
\begin{eqnarray}
\inf_{\overline{b} \in L^p(\omega;\mathbb{R}^3)} \int_{\omega}\mathcal{Q}^*W(x_{\alpha};D_{\alpha}u|\overline{b})dx_{\alpha} & \leq & \int_{\omega}\mathcal{Q}^*W(x_{\alpha};D_{\alpha}u|\overline{b}_0)dx_{\alpha}\nonumber\\
 & = & \int_{\omega}\underline{W}(x_{\alpha};D_{\alpha}u)dx_{\alpha}\nonumber\\
 & = & \int_{\omega}\min_{z \in \mathbb{R}^3}\mathcal{Q}^*W(x_{\alpha};D_{\alpha}u|z)dx_{\alpha}\nonumber\\
 & \leq & \int_{\omega}\mathcal{Q}^*W(x_{\alpha};D_{\alpha}u|\overline{b})dx_{\alpha},\label{r6}
\end{eqnarray}
where the last inequality holds for all $\overline{b} \in L^p(\omega;\mathbb{R}^3)$. Taking the infimum in $\overline{b}$ in the last term of (\ref{r6}), the inequalities become equalities thus
$$\inf_{\overline{b} \in L^p(\omega;\mathbb{R}^3)} \int_{\omega}\mathcal{Q}^*W(x_{\alpha};D_{\alpha}u|\overline{b})dx_{\alpha}=\int_{\omega}\underline{W}(x_{\alpha};D_{\alpha}u)dx_{\alpha}.$$
This shows that Theorem \ref{t2.1} is recovered from Theorem \ref{t3.1}.

\hfill October 8, 2003
\bigskip

\noindent(*, **) \small{LPMTM, Institut Galil\'ee, Universit\'e Paris-Nord,
93430 Villetaneuse, France,}

\bigskip
\noindent \small{ email(*): jfb@galilee.univ-paris13.fr\hfill email(**): francfor@galilee.univ-paris13.fr}

\addcontentsline{toc}{section}{References}

\begin{thebibliography}{50}

\bibitem{Ac&Fu} E. Acerbi, N. Fusco : \textit{Semicontinuity results in the calculus of variations}, Arch. Rat. Mech. Anal., \textbf{86}, 1984, 125-145.

\bibitem{Bo&Fo} M. Bocea, I. Fonseca : \textit{Equi-integrability results for 3D-2D dimension reduction problems}, ESAIM : Control, Optimisation and Calculus of Variations, \textbf{7}, 2002, 443-470.

\bibitem{Bo&Fo&Ma} G. Bouchitt\'e, I. Fonseca, M.L. Mascarenhas : \textit{Bending moment in membrane theory}, Preprint CNA 2002.

\bibitem{Braides} A. Braides : personal communication.

\bibitem{Bra&De} A. Braides, A. Defranceschi : Homogenization of multiple integrals, Oxford lectures series in mathematics and its applications, Clarendon Press, Oxford, 1998.

\bibitem{Bra&Fo&Fr} A. Braides, I. Fonseca, G. Francfort : \textit{3D-2D asymptotic analysis for inhomogeneous thin films}, Indiana Univ. Math. J.,\textbf{49}, 2000, 1367-1404.

\bibitem{Dac1} B. Dacorogna : Direct methods in the calculus of variations, Springer-Verlag, Berlin, 1988.

\bibitem{DalMaso} G. Dal Maso : An introduction to $\Gamma$-convergence, Birkha\"user, Boston, 1993.

\bibitem{Ek&Te} I. Ekeland, R. Temam : Analyse convexe et probl\`emes variationnels, Dunod, Gauthiers-Villars, Paris, 1974.

\bibitem{Ev&Ga} L.C. Evans, R.F. Gariepy : Measure theory and fine properties of functions, Boca Raton, CRC Press, 1992.

\bibitem{Fox&Ra&Si} D. Fox, A. Raoult, J.C. Simo : \textit{A justification of nonlinear properly invariant plate theories}, Arch. Rat. Mech. Anal., {\bf 25}, 1992, 157-199.

\bibitem{Fr&Ja&Mu1} G. Friesecke, R.D. James, S. M\"uller : \textit{Rigorous derivation of nonlinear plate theory and geometric rigidity}, C.R. Acad. Sci. Paris, S\'erie I, \textbf{334}, 2001, 173-178.

\bibitem{Fr&Ja&Mu2} G. Friesecke, R.D. James, S. M\"uller : \textit{A Theorem on geometric rigidity and the derivation of nonlinear plate theory from three dimensional elasticity},  Comm. Pure Appl. Math. 55 (2002), {\bf 11}, 1461--1506.

\bibitem{Fr&Ja&Mu3} G. Friesecke, R.D. James, S. M\"uller : \textit{The F\"oppl-von K\'arm\'an plate theory as a low energy $\Gamma$-limit of nonlinear elasticity}, C.R. Acad. Sci. Paris, S\'erie I, \textbf{335}, 2002, 201-206.

\bibitem{Le&Ra} H. Le Dret, A. Raoult : \textit{The nonlinear membrane model as variational limit of nonlinear three-dimensional elasticity}, J. Math. Pures Appl., \textbf{74}, 1995, 549-578.

\end{thebibliography}
\end{document}